%% file: article.tex
\documentclass[11pt]{article}
\usepackage[a4paper,margin=1.5cm]{geometry}
\usepackage{setspace}

\usepackage{amstext}
\usepackage{amsthm}
\usepackage{amsopn}
\usepackage{amsfonts}
\usepackage{amsmath}
\usepackage{amssymb}
\usepackage{framed}
\usepackage{enumitem}

\usepackage[dvipsnames]{xcolor}
\definecolor{greytext}{gray}{0.5}

\usepackage{url}
\usepackage[numbers,sort&compress]{natbib}

\DeclareUrlCommand\DOI{}

\usepackage{setspace}
\usepackage[bookmarks,bookmarksnumbered,linktocpage]{hyperref}

\usepackage[capitalise]{cleveref}
\usepackage[page]{appendix}
\input{cmds-article}

\newcommand{\manref}[2]{\hyperref[#1]{#2~\ref*{#1}}}

\usepackage{datetime}
\title{Homotopy theories of $(\infty, \infty)$-categories as universal fixed points with respect to enrichment}
\date{\monthname~\the\year}
\author{Zach Goldthorpe}

\begin{document}
\maketitle
\begin{abstract}
    We show that both the $\infty$-category of $(\infty, \infty)$-categories with inductively defined equivalences, and with coinductively defined equivalences, satisfy universal properties with respect to weak enrichment in the sense of Gepner and Haugseng.
    In particular, we prove that $(\infty, \infty)$-categories with coinductive equivalences form a terminal object in the $\infty$-category of fixed points for enrichment, and that $(\infty, \infty)$-categories with inductive equivalences form an initial object in the subcategory of locally presentable fixed points.
    To do so, we develop an analogue of Ad\'amek's construction of free endofunctor algebras in the $\infty$-categorical setting.
    We prove that $(\infty, \infty)$-categories with coinductive equivalences form a terminal coalgebra with respect to weak enrichment, and $(\infty, \infty)$-categories with inductive equivalences form an initial algebra with respect to weak enrichment.
\end{abstract}

\tableofcontents{}

\input{src/introduction}

\input{src/algebra}

\input{src/womcats}

\begin{appendices}
\input{src/appendix}

\end{appendices}

\addcontentsline{toc}{section}{References}
\bibliographystyle{alphaurl}
\bibliography{biblio}

\end{document}

%% file: cmds-article.tex
\usepackage[dvipsnames]{xcolor}
\definecolor{greytext}{gray}{0.5}

\usepackage{amstext}
\usepackage{amsthm}
\usepackage{amsopn}
\usepackage{amsfonts}
\usepackage{amsmath}
\usepackage{amssymb}
\usepackage{stackrel}
\usepackage{framed}

\usepackage{array}
\usepackage{dsfont}

\usepackage{mathrsfs}
\usepackage{mathtools}
\usepackage{enumitem}
\usepackage{tikz-cd}
\usepackage{arydshln}
\usepackage{subcaption}

\usepackage[scr = euler]{mathalpha}

\usepackage{chngcntr}
\counterwithin*{equation}{section}
\counterwithin*{figure}{section}

\newtheorem{maintheorem}{Theorem}

\newtheorem{theorem}{Theorem}[subsection]
\newtheorem{corollary}[theorem]{Corollary}
\newtheorem{lemma}[theorem]{Lemma}
\newtheorem{proposition}[theorem]{Proposition}

\theoremstyle{definition}
\newtheorem{definition}[theorem]{Definition}
\newtheorem{example}[theorem]{Example}

\theoremstyle{remark}
\newtheorem{remark}[theorem]{Remark}

\let\oldcite\cite
\renewcommand\cite[2][]{{\rm\oldcite[#1]{#2}}}

\numberwithin{equation}{section}

\def\(#1\){$\displaystyle#1$}
\def\[#1\]{\begin{align*}#1\end{align*}}

\renewcommand\bar[1]{\overline{#1}}
\renewcommand\hat[1]{\widehat{#1}}

\newcommand\1{\mathds{1}}

\newcommand\E{\mathbb{E}}

\newcommand\bI{\mathbb{I}}

\newcommand\bS{\mathbb{S}}
\newcommand\Z{\mathbb{Z}}

\newcommand\bfB{\mathbf{B}}

\newcommand\bGamma{\mathbf{\Gamma}}
\newcommand\bDelta{\mathbf{\Delta}}

\newcommand\bTheta{\mathbf{\Theta}}

\newcommand\cC{\mathcal{C}}
\newcommand\cD{\mathcal{D}}

\newcommand\cV{\mathcal{V}}

\newcommand\sA{\mathscr{A}}
\newcommand\sB{\mathscr{B}}
\newcommand\sC{\mathscr{C}}
\newcommand\sD{\mathscr{D}}

\newcommand\sK{\mathscr{K}}
\newcommand\sL{\mathscr{L}}

\newcommand\sO{\mathscr{O}}

\newcommand\sS{\mathscr{S}}
\newcommand\sV{\mathscr{V}}

\newcommand\frakM{\mathfrak{M}}
\newcommand\frakO{\mathfrak{O}}
\newcommand\frakP{\mathfrak{P}}

\newcommand\cat{\mathrm{cat}}

\newcommand\co{\mathrm{co}}

\newcommand\gen{\mathrm{gen}}

\newcommand\id{\operatorname{id}}
\newcommand\lax{\mathrm{lax}}

\newcommand\op{\mathrm{op}}

\newcommand\rank{\operatorname{rank}}

\newcommand\sk{\operatorname{sk}}

\newcommand\Hom{\operatorname{Hom}}
\newcommand\Id{\operatorname{Id}}

\newcommand\Map{\operatorname{Map}}
\newcommand\Nat{\operatorname{Nat}}
\newcommand\Ob{\operatorname{Ob}}

\newcommand\varlim{\varprojlim}
\newcommand\varcolim{\varinjlim}

\makeatletter
\newcommand{\pto}{}
\newcommand{\pgets}{}

\DeclareRobustCommand{\pto}{\mathrel{\mathpalette\p@to@gets\to}}
\DeclareRobustCommand{\pgets}{\mathrel{\mathpalette\p@to@gets\gets}}

\newcommand{\p@to@gets}[2]{
  \ooalign{\hidewidth$\m@th#1\mapstochar\mkern5mu$\hidewidth\cr$\m@th#1\to$\cr}
}
\makeatother

\newcommand\Alg{\mathbf{Alg}}

\newcommand\Cat{\mathbf{Cat}}

\newcommand\Enr{\mathbf{Enr}}

\newcommand\Fix{\mathbf{Fix}}
\newcommand\Fun{\mathbf{Fun}}
\newcommand\Grpd{\mathbf{Grpd}}

\newcommand\Mon{\mathbf{Mon}}
\newcommand\Opd{\mathbf{Opd}}
\newcommand\Pres{\mathbf{Pres}}

\newcommand\Set{\mathbf{Set}}
\newcommand\sSet{\mathbf{sSet}}

\newcommand\SymMon{\mathbf{SymMon}}

%% file: src/introduction.tex
%!TEX root= ../article.tex

\section{Introduction}\label{sec:intro}

Following the extensive development of the theory of quasicategories, especially by Joyal and Lurie, the study of many mathematical objects in the language of higher category theory has become much more convenient.
Quasicategories serve as a homotopy coherent generalisation of ordinary categories, and the resulting theory can, in many cases, be viewed in precisely this way.
This allows a practising mathematician to work with these models for higher category theory completely analogously to how they would have worked with ordinary categories, but with the advantage that they can study objects with more nuanced (i.e., homotopical or derived) structure.

Higher category theory is a generalisation of ordinary category theory, where morphisms can be compared via other morphisms of higher dimension.
Specifically, a higher category $\sC$ would consist of a collection of objects and a collection of morphisms between objects, as in an ordinary category, but also include $2$-morphisms between morphisms, and $3$-morphisms between $2$-morphisms, \emph{ad infinitum}.

In particular, between any two objects $x, y$ of a higher category $\sC$ is a higher category $\Hom_\sC(x,y)$ of morphisms between them, where the objects are morphisms $x\to y$, the morphisms are $2$-morphisms in $\sC$, and so on.
This perspective reveals the relationship between higher category theory and another generalisation of category theory: \emph{enriched} category theory.
Briefly, a category $\cC$ is said to be \emph{enriched} in another category $\cV$ if we can realise the collection of morphisms between any pair $x, y$ of objects of $\cC$ as an object $\Hom_\cC(x, y)$ of $\cV$.
The classical theory of categories enriched in ordinary categories was developed extensively by Kelly in \cite{kelly}.

By enriching a category in a higher category $\sV$, we should then expect a higher category with morphisms of dimension up to one greater than the highest dimension of morphism in $\sV$.
More specifically, say that a higher category $\sC$ is an \emph{$(n,r)$-category} if morphisms of dimension at least $n+1$ are trivial, and morphisms of dimension at least $r+1$ are invertible.
For instance, $(1, 1)$-categories are ordinary categories, and $(1, 0)$-categories are ordinary groupoids.
The expectation is then that all $(n+1, r+1)$-categories can be realised via enrichment in $(n, r)$-categories.

However, the classical theory of enrichment fails to encompass the desired flexibility of higher category theory.
The composition of morphisms in a category enriched in an ordinary category is necessarily strictly associative, so higher category theory defined through classical enrichment would likewise have strictly associative composition between morphisms.
Such higher categories are referred to as \emph{strict} higher categories, and turn out to be too rigid for a suitable theory of higher categories in application.

Therefore, a suitable approach to higher category theory via enrichment requires a theory of enrichment capable of expressing associativity of composition up to an enriched analogue of coherent homotopy.
This is difficult due to the circular nature of the construction: the most suitable context for coherent homotopy is inside of a higher category, and thus a fully general theory of enrichment is impossible without a suitably general theory of higher category theory.
Nonetheless, this enrichment paradigm in tandem with the Homotopy Hypothesis leads to a means for modelling subclasses of higher categories: if we accept the Homotopy Hypothesis, then we can take the homotopy types of spaces as a model for higher groupoids or $(\infty, 0)$-categories.
If one can make sense of weak enrichment in $(\infty, 0)$-categories, then one obtains a theory of $(\infty, 1)$-categories.

This philosophy leads to several of the popular models of $(\infty, 1)$-categories, among which include complete Segal spaces, Segal categories, and quasicategories.
Thanks to the work of Joyal and Lurie (largely compiled into \cite{lurie}), quasicategories have become one of the most popular models of $(\infty, 1)$-category theory, but many of the theories of $(\infty, 1)$-categories have been shown to be equivalent.

Exploration of enrichment in these models for $(\infty, 1)$-categories leads to models for $(\infty, n)$-categories for $n\geq1$ such as complete $n$-fold Segal spaces, complete $\bTheta_n$-spaces, Segal $n$-categories (all described and proven equivalent in \cite{barwick-schommer-pries}), and $n$-quasicategories (developed and proven equivalent to $\bTheta_n$-spaces in \cite{ara}), all of which provide equivalent models for the $(\infty,1)$-category of $(\infty, n)$-categories.
In fact, the Unicity Theorem \cite[Theorem 7.3]{barwick-schommer-pries} proves that there is an essentially unique reasonable $(\infty,1)$-category of $(\infty, n)$-categories, and the moduli space of such theories is a $B(\Z/2)^n$ for all $0\leq n<\infty$.

Although not the most general class of higher categories, $(\infty, 1)$-categories are a natural context for studying categorical constructions in tandem with abstract homotopy theory, making them a suitably general context for many applications of higher category theory to mathematics.
In this context, Gepner and Haugseng develop in \cite{gepner-haugseng} a fully general theory of enrichment in a monoidal $(\infty, 1)$-category $\sV$, wherein composition in a category enriched in $\sV$ is only associative up to coherent homotopy in $\sV$.
In particular, this general framework of enrichment provides a uniform construction of $(\infty, n)$-categories for all finite $n$, which is shown to be equivalent to $n$-fold Segal spaces in \cite{haugseng}.

Up to a choice of suitable model at every dimension, we have a tower of inclusions
\[
    \Cat_{(\infty,0)} \hookrightarrow \Cat_{(\infty,1)} \hookrightarrow \Cat_{(\infty,2)} \hookrightarrow \Cat_{(\infty,3)} \hookrightarrow \dots
\]
Each inclusion $\Cat_{(\infty, n)}\hookrightarrow\Cat_{(\infty,n+1)}$ admits both a left adjoint $\pi$ and a right adjoint $\kappa$.
The left adjoint acts on an $(\infty,n+1)$-category $\sC$ by producing the $(\infty, n)$-category $\pi\sC$ obtained by formally inverting all of the $(n+1)$-morphisms, whereas the right adjoint acts on $\sC$ by producing the maximal sub-$(\infty,n)$-category $\kappa\sC$ of $\sC$.
These adjoints provide constructions of two candidates for the $(\infty,1)$-category of fully weak higher categories:
\[
    \Cat_\omega &:= \varlim\left(\dots\to\Cat_{(\infty,3)}\xrightarrow\pi\Cat_{(\infty,2)}\xrightarrow\pi\Cat_{(\infty,1)}\xrightarrow\pi\Cat_{(\infty,0)}\right) \\
    \Cat_{(\infty, \infty)} &:= \varlim\left(\dots\to\Cat_{(\infty,3)}\xrightarrow\kappa\Cat_{(\infty,2)}\xrightarrow\kappa\Cat_{(\infty,1)}\xrightarrow\kappa\Cat_{(\infty,0)}\right)
\]
Roughly speaking, the difference between the $(\infty,1)$-category $\Cat_\omega$ of \emph{$\omega$-categories} and the $(\infty,1)$-category $\Cat_{(\infty,\infty)}$ of \emph{$(\infty, \infty)$-categories} is that the equivalences in an $\omega$-category are defined weakly coinductively, whereas equivalences in an $(\infty, \infty)$-category are defined inductively from the higher identity endomorphisms.
The different notions of equivalence are described in more detail in \cref{rem:ciwe}.

In this paper, we show that the latter theory $\Cat_{(\infty, \infty)}$ is compatible with enrichment.
%In this paper, we show that both theories are compatible with enrichment.
For a symmetric monoidal $(\infty, 1)$-category $\sV$, let $\sV\Cat$ denote the symmetric monoidal $(\infty, 1)$-category of categories weakly enriched in $\sV$.
Then, we show that there is a canonical equivalence $\Cat_{(\infty, \infty)} \simeq (\Cat_{(\infty, \infty)})\Cat$.
%Then, we show that there are canonical equivalences $\Cat_\omega\simeq(\Cat_\omega)\Cat$ and $\Cat_{(\infty,\infty)}\simeq(\Cat_{(\infty,\infty)})\Cat$.
In fact, we prove a stronger result:
\begin{maintheorem}[\cref{thm:main}]\label{thm:A}
    The $(\infty, 1)$-category %$\Cat_\omega$ defines a terminal object in the $(\infty, 1)$-category of symmetric monoidal $(\infty,1)$-categories $\sV$ with a symmetric monoidal equivalence $\sV\xrightarrow\sim\sV\Cat$.
    %On the other hand,
    $\Cat_{(\infty,\infty)}$ defines an initial object in the $(\infty, 1)$-category of presentably symmetric monoidal $(\infty, 1)$-categories $\sV$ with a symmetric monoidal equivalence $\sV\xrightarrow\sim\sV\Cat$.
\end{maintheorem}

In order to prove this result, we study (co)algebras associated to any endofunctor of a general $(\infty,1)$-category, and prove an $(\infty, 1)$-categorical analogue of Ad\'amek's construction of an initial endofunctor algebra (see \cite[Theorem 3.19]{adamek:modern}).
More precisely, fix an $(\infty, 1)$-category $\sK$ and a functor $F:\sK\to\sK$, then an \emph{$F$-algebra} is a pair $(A,\alpha)$, where $A$ is an object of $\sK$, and $\alpha:FA\to A$ is a morphism of $\sK$ called the \emph{action}.
The $F$-algebras in $\sK$ collect to form an $(\infty, 1)$-category $\sK(F)$.
Then, we prove:

\begin{maintheorem}[\cref{thm:adamek}]\label{thm:B}
    Let $\emptyset$ be an initial object of $\sK$.
    If the diagram
    \[
        \emptyset \xrightarrow{!}F\emptyset \xrightarrow{F(!)} F^2\emptyset \xrightarrow{F^2(!)} F^3\emptyset \to \dots
    \]
    admits a colimit $I$, and the canonical map $f:I\to FI$ is an equivalence, then $(I,f^{-1})$ is an initial $F$-algebra.
\end{maintheorem}

Enrichment restricts to an endofunctor on (presentably) symmetric monoidal $(\infty, 1)$-categories, so \cref{thm:A} ultimately follows from %observing that $\Cat_\omega$ defines a terminal coalgebra for enrichment over symmetric monoidal $(\infty,1)$-categories, and $\Cat_{(\infty,\infty)}$ defines an initial algebra for enrichment over presentably symmetric monoidal $(\infty, 1)$-categories.
demonstrating that $\Cat_{(\infty,\infty)}$ defines an initial algebra for enrichment over presentably symmetric monoidal $(\infty, 1)$-categories.

\input{src/introduction/organisation}

\input{src/introduction/notation}

%% file: src/introduction/organisation.tex
%%!TEX root = ../../article.tex

\subsection{Organisation of paper}\label{ssec:org}

In \cref{sec:alg}, we develop an analogue of the theory of (co)algebras of an endofunctor in the context of $(\infty, 1)$-categories.
In order to describe Ad\'amek's construction and prove its correctness, we find it is more convenient to embed endofunctor algebras into a larger category of \emph{lax} algebras, which we define in \cref{ssec:defs}.
In \cref{ssec:adamek}, we generalise Ad\'amek's construction to produce free endofunctor algebras generated by a general lax algebra.
By specialising this construction to lax algebras freely generated by a single object, we recover the usual construction in an $(\infty,1)$-categorical context, and in particular we prove \cref{thm:B}.
In \cref{ssec:fix}, we use the developed theory of endofunctor (co)algebras to study general fixed points of an endofunctor.
In particular, we show in \cref{thm:fixloc} how to freely construct a fixed point when given a coalgebra for any endofunctor satisfying relatively flexible constraints (for instance, it suffices if the endofunctor is accessible and lives over a locally presentable $(\infty, 1)$-category).

In \cref{ssec:operad,ssec:endo}, we review the necessary groundwork of \cite{gepner-haugseng} for studying enrichment in a monoidal $(\infty, 1)$-category.
We recall that enrichment restricts to define an endofunctor on (presentably) symmetric monoidal $(\infty, 1)$-categories, and prove in \cref{thm:wclim} that enrichment preserves suitably well-behaved limits of diagrams of monoidal $(\infty, 1)$-categories indexed by weakly contractible simplicial sets.
This puts us in a position to employ the constructions of \cref{ssec:fix}, which we do in \cref{ssec:womcat} to prove our main result \cref{thm:A}.

%% file: src/introduction/notation.tex
%%!TEX root = ../../article.tex

\subsection{Notation and terminology}\label{ssec:not}

Let $\bDelta$ denote the simplex category, and let $\sSet := \Fun(\bDelta^\op,\Set)$ be the category of simplicial sets.
For a simplicial set $S:\bDelta^\op\to\Set$ and $[n]$ an object of $\bDelta$, let $S_n := S([n])$ denote the set of $n$-cells of $S$.

Denote by $\Delta[n] : \bDelta^\op\to\Set$ the functor represented by the object $[n]$ in $\bDelta$.
Then, define $\partial\Delta[n]$ to be the simplicial subset of $\Delta[n]$ obtained by removing the unique degenerate $n$-cell, and by $\Lambda^i[n]$ for $0\leq i\leq n$ the simplicial subset of $\partial\Delta[n]$ obtained by removing the $(n-1)$-cell associated to the inclusion $[n-1]\hookrightarrow[n]$ whose image does not contain $i$.

Call a simplicial set $S$ a \emph{space} (or an $\infty$-\emph{groupoid}) if it is a Kan complex; that is, $S$ satisfies the right lifting property against all horn inclusions $\Lambda^i[n]\hookrightarrow\Delta[n]$ for $0\leq i\leq n$.
Call a simplicial set $\sC$ an $\infty$-\emph{category} (which in this paper is synonymous with an $(\infty, 1)$-category) if it is a weak Kan complex; that is, $\sC$ satisfies the right lifting property against all inner horn inclusions $\Lambda^i[n]\hookrightarrow\Delta[n]$ for $0<i<n$.

We denote by $\sS$ the large $\infty$-category of spaces; that is, the $\infty$-category induced by Quillen's model structure on $\sSet$.
When we want to view these spaces as models for $\infty$-groupoids, we may also denote this category by $\Grpd_\infty := \sS$.
Denote by $\Cat_\infty$ the large $\infty$-category of $\infty$-categories, induced by the Joyal model structure on $\sSet$.

In an $\infty$-category $\sC$, we just refer to $\infty$-limits and $\infty$-colimits as limits and colimits, respectively.
Call a functor $F:\sC\to\sD$ of $\infty$-categories \emph{continuous} if it preserves limits, and \emph{cocontinuous} if it preserves colimits.
$F$ is said to \emph{reflect limits} if whenever $\bar p:K^\triangleleft\to\sC$ is a diagram such that $F\circ\bar p$ is a limit diagram, then $\bar p$ is a limit diagram.
If $F$ preserves and reflects limits, then $F$ is said to \emph{create limits}.
The notions of reflecting and creating colimits are entirely analogous.

%% file: src/algebra.tex
%%!TEX root = ../article.tex

\section{Algebras of an \texorpdfstring{$\infty$}{oo}-endofunctor}\label{sec:alg}

Fix an $\infty$-category $\sK$ and an endofunctor $F:\sK\to\sK$.
We are primarily interested in studying the fixed points of $F$; that is, the objects of $\sK$ that are equivalent to their image under $F$.
This is intimately related to the theory of (co)algebras of $F$ studied (in the classical setting) by Lambek in \cite{lambek} and Ad\'amek in \cite{adamek, adamek:modern}.

An $F$-\emph{algebra} is a pair $(A,\alpha)$, where $A$ is an object of $\sK$, and $\alpha:FA\to A$ is a morphism called the ``action'' of $A$.
An $F$-\emph{algebra homomorphism} is then a morphism of underlying objects that commutes with the $F$-algebra actions up to homotopy.

\begin{definition}\label{def:K(F)}
    Define the $\infty$-category $\sK(F)$ of $F$-algebras and $F$-algebra homomorphisms as the pullback
    $$
    \begin{tikzcd}
        \sK(F) \ar[r]\ar[d]\ar[dr, phantom, "\lrcorner" very near start]
        & \sK^{\Delta[1]} \ar[d, two heads] \\
        \sK \ar[r, "{(F,\id)}"']
        & \sK\times\sK
    \end{tikzcd}
    $$
    Note that the vertical map on the right is induced by the inclusion $\{0,1\}\hookrightarrow\Delta[1]$, and is thus a categorical fibration.
    Therefore, the pullback is a homotopy pullback of $\infty$-categories.
\end{definition}

We can then identify the fixed points of $F$ with the $F$-algebras with invertible actions.
Lambek observed in \cite[Lemma 2.2]{lambek} that initial $F$-algebras are always fixed points of $F$:
\begin{lemma}[Lambek]\label{lem:lambek}
    Suppose $(I,i)$ is an initial object in the $\infty$-category $\sK(F)$.
    Then, the action $i:FI\to I$ is an equivalence.
\end{lemma}
\begin{proof}
    Consider the $F$-algebra $(FI,Fi)$.
    Since $(I,i)$ is initial, there is an essentially unique $F$-algebra homomorphism $u:(I,i)\to(FI,Fi)$.
    The composite $i\circ u$ thus defines an $F$-algebra endomorphism of $(I,i)$, and must therefore be homotopic to the identity; that is, $i\circ u\simeq\id_I$.

    On the other hand, consider the diagram
    $$
    \begin{tikzcd}[sep=large]
        FI \ar[r, "Fu"]\ar[d, "i"']\ar[dr, "{F(i\circ u)}"]
        & FFI \ar[d, "Fi"] \\
        I \ar[r, "u"']
        & FI
    \end{tikzcd}
    $$
    The perimeter commutes up to homotopy because $u$ is an $F$-algebra homomorphism, and the upper triangle commutes by the functoriality of $F$.
    Since $i\circ u\simeq\id_I$, it follows that also $u\circ i\simeq F(i\circ u)\simeq \id_{FI}$.
    Therefore, $u\simeq i^{-1}$, proving that $i$ is an equivalence, as desired.
\end{proof}

In \cref{ssec:adamek}, we prove an $\infty$-categorical analogue of Ad\'amek's construction of free $F$-algebras generated by objects of $\sK$ (see \cite[p. 592]{adamek}).
In particular, by taking the free $F$-algebra generated by an initial object of $\sK$, we produce an initial $F$-algebra, and thus a universal fixed point by \hyperref[lem:lambek]{Lambek's lemma} above.

Note that the entire theory of $F$-algebras dualises to give a theory of $F$-\emph{coalgebras}.
Explicitly, an $F$-coalgebra is a pair $(C,\nu)$, where $C$ is an object of $\sK$, and $\nu:C\to FC$ is a morphism called the ``coaction,'' and an $F$-coalgebra homomorphism is a morphism of the underlying objects that commutes with the coactions.
If we denote the $\infty$-category of $F$-coalgebras and $F$-coalgebra homomorphisms by $\sK_\co(F)$, then we have an equivalence $\sK_\co(F)\simeq(\sK^\op(F))^\op$.
Therefore, the theories of $F$-algebras and $F$-coalgebras are entirely dual.

We will primarily develop the theory for $F$-algebras, but $F$-coalgebras will play a more prominent role when studying general fixed points in \cref{ssec:fix}.

\input{src/algebra/definitions}

\input{src/algebra/adamek}

\input{src/algebra/fixedpoints}

%% file: src/algebra/definitions.tex
%%!TEX root = ../../article.tex

\subsection{Algebras and lax algebras}\label{ssec:defs}

Ad\'amek's construction of a free $F$-algebra defines the underlying object as a colimit, and the action is given as the inverse of a canonical map induced by the universal property of this colimit.
In the context of $\infty$-categories, an inverse is only unique up to homotopy, which makes the classical proof difficult to replicate in this setting.
Therefore, we embed $\sK(F)$ into a larger category $\sK^\lax(F)$ of \emph{lax $F$-algebras}.

\begin{definition}\label{def:Klax}
    Define the $\infty$-category $\sK^\lax(F)$ as the pullback
    $$
    \begin{tikzcd}
        \sK^\lax(F) \ar[r]\ar[d]\ar[dr, phantom, "\lrcorner" very near start]
        & \sK^{\{1\gets0\to2\}} \ar[d, two heads] \\
        \sK \ar[r, "{(F,\id)}"']
        & \sK^{\{1\}}\times\sK^{\{2\}}
    \end{tikzcd}
    $$
    The vertical map on the right is induced by the cofibration $\{1,2\}\hookrightarrow\{1\gets0\to2\}$, and is thus a categorical fibration, showing that the pullback square is a homotopy pullback of $\infty$-categories.
\end{definition}
\begin{remark}\label{rem:tower}
    Every square in the tower below is a pullback square:
    $$
    \begin{tikzcd}
        \sK^\lax(F) \ar[r]\ar[d]\ar[dr, phantom, "\lrcorner" very near start]
        & \sK^{\{1\gets0\to2\}} \ar[d, two heads] \\
        \sK^{\{0\to2\}} \ar[r, "{(F(2),\id)}"]\ar[d]\ar[dr, phantom, "\lrcorner" very near start]
        & \sK^{\{1\}}\times\sK^{\{0\to2\}} \ar[d, two heads] \\
        \sK^{\{0\}}\times\sK^{\{2\}} \ar[r, "{(F(2),\id)}"]\ar[d]\ar[dr, phantom, "\lrcorner" very near start]
        & \sK^{\{1\}}\times\sK^{\{0\}}\times\sK^{\{2\}} \ar[d] \\
        \sK^{\{2\}} \ar[r, "{(F,\id)}"']
        & \sK^{\{1\}}\times\sK^{\{2\}}
    \end{tikzcd}
    $$
    As the vertical maps on the right are induced by inclusions of simplicial sets, it follows that they are categorical fibrations, showing that all of these pullback squares are also homotopy pullback squares of $\infty$-categories.
\end{remark}

Concretely, a lax $F$-algebra is a span $FB\xleftarrow rE\xrightarrow aB$ in $\sK$, where the morphism $a$ may be called the ``lax action,'' and the morphism $r$ may be called a ``resolution.''

\begin{proposition}\label{prop:embed}
    Let $\hat\sK(F)$ denote the full subcategory of $\sK^\lax(F)$ spanned by those lax $F$-algebras $FB\gets E\to B$ where the resolution $E\to FB$ is invertible.
    Then, there is a canonical equivalence of $\infty$-categories $\sK(F)\xrightarrow\sim\hat\sK(F)$.
\end{proposition}
\begin{proof}
    Note that $\hat\sK(F)$ can be defined equivalently as follows, using the cartesian model structure on marked simplicial sets of \cite[\S3.1]{lurie}.
    Let $\Lambda_+^0[2] := \{1\xleftarrow+0\to2\}$ denote the marked simplicial set obtained by taking the walking span $\Lambda^0[2]$ and marking the left-pointing edge.
    Then, $\hat\sK(F)$ is the pullback
    $$
    \begin{tikzcd}
        \hat\sK(F) \ar[r]\ar[d]\ar[dr, phantom, "\lrcorner" very near start]
        & \Map^\flat(\Lambda^0_+[2],\sK^\natural) \ar[d, two heads] \\
        \sK \ar[r, "{(F,\id)}"']
        & \sK\times\sK
    \end{tikzcd}
    $$
    where $\Map^\flat(X,Y)$ is the underlying simplicial set of the internal hom of marked simplicial sets, and $\sK^\natural$ is the $\infty$-category $\sK$ marked at the equivalences.
    By \cite[Remark 3.1.4.5]{lurie}, $\Map^\flat$ provides the cartesian model structure with an enrichment in Joyal's model structure on simplicial sets.
    Note that $\sK^\natural$ is fibrant in the cartesian model structure, and so the vertical map induced by an inclusion of marked simplicial sets is therefore a categorical fibration.

    Let $\Delta[1]^\flat$ denote the simplicial set $\Delta[1]$ marked only at the degenerate edges.
    Then, the inclusion $\Delta[1]^\flat\to\Lambda_+^0[2]$ that picks out the unmarked edge of $\Lambda_+^0[2]$ is marked anodyne by \cite[Proposition 3.1.1.5]{lurie}, and admits a retraction $\Lambda_+^0[2]\to\Delta[1]$.
    By 2-out-of-3, it follows that this retraction is a cartesian equivalence.
    Therefore, since $\sK^\natural$ is fibrant, the retraction induces a categorical equivalence $\sK^{\Delta[1]}=\Map^\flat(\Delta[1]^\flat,\sK^\natural)\to\Map^\flat(\Lambda_+^0[2],\sK^\natural)$.
    In particular, we have the following diagram:
    $$
    \begin{tikzcd}
        \sK(F) \ar[rr]\ar[ddr, bend right]\ar[dr, dashed]
        && \sK^{\Delta[1]} \ar[d, "\sim" sloped] \\
        & \hat\sK(F) \ar[r]\ar[d]\ar[dr, phantom, "\lrcorner" very near start]
        & \Map^\flat(\Lambda^0_+[2],\sK^\natural) \ar[d, two heads] \\
        & \sK \ar[r, "{(F,\id)}"']
        & \sK\times\sK
    \end{tikzcd}
    $$
    The perimeter is the definitional pullback square for $\sK(F)$, so both the perimeter and the inner pullback square are homotopy pullbacks.
    Since the objects of the two pullback diagrams are connected by categorical equivalences, it thus follows that the induced map $\sK(F)\to\hat\sK(F)$ is a categorical equivalence as well.
\end{proof}

The embedding identifies an $F$-algebra $(A,\alpha)$ with the lax $F$-algebra $FA=FA\xrightarrow\alpha A$, and conversely any lax $F$-algebra $FB\xleftarrow rE\xrightarrow aB$ with an invertible resolution induces an $F$-algerba $(B,ar^{-1})$.
We therefore tacitly identify $\sK(F)$ with its essential image in $\sK^\lax(F)$.

In the context of \hyperref[thm:adamek]{Ad\'amek's construction}, we can avoid explicitly inverting the canonical map $I\to FI$ to define the action, and instead prove that the lax $F$-algebra $FI\xleftarrow\sim I = I$ is initial in the essential image of $\sK(F)$.

%% file: src/algebra/adamek.tex
%%!TEX root = ../../article.tex

\subsection{\texorpdfstring{Ad\'amek's}{Adamek's} construction}\label{ssec:adamek}

In this section, we generalise Ad\'amek's construction to produce a free $F$-algebra generated by any lax $F$-algebra $FB\gets E\to B$.
If $\sK$ has an initial object $\emptyset$, then the forgetful functor $\sK^\lax(F)\to\sK$ sending a lax $F$-algebra $FB\gets E\to B$ to the underlying object $B$ admits a left adjoint given by sending an object $K$ of $\sK$ to the lax $F$-algebra $FK\gets\emptyset\to K$.
Therefore, we recover an $\infty$-categorical analogue of Ad\'amek's construction of a free $F$-algebra generated by an object $K$ by applying our generalised construction to the lax $F$-algebra $FK\gets\emptyset\to K$.

Classically, Ad\'amek's free $F$-algebra construction on an object $K$ is given by taking the colimit of a diagram
\[
    K\to K\sqcup FK \to K\sqcup F(K\sqcup FK) \to K\sqcup F(K\sqcup F(K\sqcup FK))
\]
The stages of this construction are computed inductively via propagation of a certain operation $X\mapsto K\sqcup FX$.
Ad\'amek shows on \cite[p. 592]{adamek} that if the colimit is preserved by this operation, then the colimit canonically carries the structure of a free $F$-algebra on $K$.

The structure of this propagation is more evident when presented as an infinite row of pushout squares:
$$
\begin{tikzcd}
    \emptyset \ar[r, "!"]\ar[d]\ar[dr, phantom, "\ulcorner" very near end]
    & FK \ar[r, "Fi_1"]\ar[d]\ar[dr, phantom, "\ulcorner" very near end]
    & F(K\sqcup FK) \ar[r, "Fi_2"]\ar[d]\ar[dr, phantom, "\ulcorner" very near end]
    & F(K\sqcup F(K\sqcup FK)) \ar[r, "Fi_3"]\ar[d]
    & \dots \\
    K \ar[r, "i_1"']
    & K\sqcup FK \ar[r, "i_2"']
    & K\sqcup F(K\sqcup FK) \ar[r, "i_3"']
    & K\sqcup F(K\sqcup F(K\sqcup FK)) \ar[r, "i_4"']
    & \dots
\end{tikzcd}
$$
Note that the first stage is precisely the pushout of the free lax $F$-algebra generated by $K$.
This suggests a natural generalisation of this ``propagation'' construction for a general lax $F$-algebra.

\begin{definition}\label{def:prop}
    Suppose $\sK$ has all finite colimits.
    For any lax $F$-algebra $FB\xleftarrow rE\xrightarrow aB$, consider the following diagram:
    $$
    \begin{tikzcd}
        FB \ar[r, "Fi"]
        & F(B\sqcup_EFB) \\
        E \ar[u, "r"]\ar[r, "r"]\ar[d, "a"']\ar[dr, phantom, "\ulcorner" very near end]
        & FB \ar[u, "Fi"']\ar[d] \\
        B \ar[r, "i"']
        & B\sqcup_EFB
    \end{tikzcd}
    $$
    The vertical arrows on the right define a new lax $F$-algebra, which we denote by $\Pi(FB\gets E\to B)$.
    This construction extends to an endofunctor $\Pi:\sK^\lax(F)\to\sK^\lax(F)$, and the horizontal arrows define a canonical natural transformation $\eta:\Id\Rightarrow\Pi$.

    We refer to $\Pi$ as the \emph{propagation} of lax $F$-algebras, and $\eta$ as the \emph{unit} of the endofunctor.
    The propagation and its unit are defined more carefully in \cref{adx:prop}.
\end{definition}

Through the propagation endofunctor and unit, Ad\'amek's construction can be seen as a special case of a more general result:

\begin{theorem}[Free fixed point construction]\label{thm:freefix}
    Let $\sL$ be an $\infty$-category, and $\Pi:\sL\to\sL$ an endofunctor with a unit; that is, with a natural transformation $\eta:\Id\Rightarrow\sL$.
    Denote by $\sL^\Pi$ the full subcategory spanned by objects $K$ such that $\eta_K:K\to\Pi K$ is an equivalence.
    
    For an ordinal $\theta$, let $[\theta]$ denote the nerve of the poset of all ordinals $0\leq\xi\leq\theta$ viewed as a category.
    For any object $L\in\sL$, construct the diagrams $D_L^\theta : [\theta]\to\sL$ by transfinite induction as follows:
    \begin{itemize}
        \item Define $D_L^0 : [0]\to\sL$ to be the diagram picking out the object $L$.
        \item Given $D_L^\theta$, define $D_L^{\theta+1}$ to be the extension of $D_L^\theta$ that sends the morphism $\theta\leq\theta+1$ in $[\theta+1]$ to $\eta_{D_L^\theta(\theta)} : D_L^\theta(\theta)\to\Pi D_L^\theta(\theta)$.
        \item For a limit ordinal $\lambda$, given $D_L^\theta$ for all $\theta<\lambda$, let $\underline\lambda = \varcolim_{\theta<\lambda}[\theta]$ be the nerve of the poset of all ordinals $0\leq\xi<\lambda$, which induces a functor $D_L^{<\lambda} : \underline\lambda \to \sL$.
            Since $[\lambda]\cong\underline\lambda^\triangleright$, define $D_L^\lambda$ to be a colimit cocone for $D_L^{<\lambda}$, if it exists.
    \end{itemize}
    Suppose for some limit ordinal $\lambda$ that the diagram $D_L^\lambda : [\lambda] \to \sL$ is well-defined, and let $\hat L := D_L^\lambda(\lambda)$.
    Then, the following are equivalent:
    \begin{enumerate}[label={(\roman*)}]
        \item\label{it:fix} $\hat L\in\sL^\Pi$.
        \item\label{it:rep} $\hat L$ corepresents the functor $\Hom_\sL(L,-)|:\sL^\Pi\to\sS$; that is, $\hat L$ is the free object in $\sL^\Pi$ generated by $L$.
    \end{enumerate}
\end{theorem}
\begin{proof}
    Since \ref{it:rep} certainly implies \ref{it:fix}, we need to show that \ref{it:fix} implies \ref{it:rep}, for which it is enough to prove that the coprojection $L\to\hat L$ induces a homotopy equivalence $\Hom_\sL(\hat L,K)\to\Hom_\sL(L,K)$ whenever $K\in\sL^\Pi$.
    Indeed, if $\hat L\in\sL^\Pi$, this would prove that $\Hom_{\sL^\Pi}(\hat L,K) = \Hom_\sL(L,K)$ for all $K\in\sL^\Pi$.
    
    Fix $K\in\sL^\Pi$ and let $\delta_K^\theta : [\theta] \to \sL$ denote the constant diagram on $K$.
    By transfinite induction, we can define a natural transformation $\delta_K^\theta\Rightarrow D_K^\theta$ for every ordinal $\theta$, where the component $K=\delta_K^\theta(\xi)\to D_K^\theta(\xi)$ is given by the transfinite composite $K\xrightarrow\eta\Pi K\xrightarrow{\Pi\eta}\Pi^2K\to\dots\to D_K^\theta(\xi)$.
    Since $K\in\sL^\Pi$, the map $\eta_K:K\to\Pi K$ is an equivalence, which ensures that the diagrams $D_K^\theta$ are well-defined for all ordinals $\theta$, and moreover that the natural transformation $\delta_K^\theta\Rightarrow D_K^\theta$ is a natural equivalence.
    
    Now, consider the diagram
    $$
    \begin{tikzcd}[ampersand replacement=\&]
        \Hom_\sL(\hat L,K) \ar[r]\ar[d, "\sim"' sloped]
        \& \Hom_\sL(L,K) \ar[d, hook] \\
        \Nat(D_L^\lambda,\delta_K^\lambda) \ar[r, "\sim"']
        \& \Nat(D_L^\lambda,D_K^\lambda) \ar[u, bend right, dashed]
    \end{tikzcd}
    $$
    The vertical map on the left is an equivalence by \cite[Lemma 4.2.4.3(ii)]{lurie}, and the horizontal map on the bottom is an equivalence since $\delta_K^\lambda\Rightarrow D_K^\lambda$ is a natural equivalence.
    The vertical map on the right is given by the functoriality of the construction of the diagram $D_{(-)}^\lambda$, and admits a retraction (denoted by the dashed arrow) that acts by projecting a natural transformation $D_L^\lambda\Rightarrow D_K^\lambda$ to the zeroth component $L\to K$.
    Since the square commutes, it follows from the 2-out-of-6 property that all of the arrows in the diagram are equivalences.
    In particular, the map $\Hom_\sL(\hat L,K)\to\Hom_\sL(L,K)$ is a weak equivalence, as desired.
\end{proof}

\begin{lemma}\label{lem:pfp}
    Let $\sK$ be finitely cocomplete so that we have the propagation endofunctor and unit on $\sK^\lax(F)$.
    Then, the inclusion $\sK(F)\hookrightarrow\sK^\lax(F)$ factors through the full subcategory $\sK^\lax(F)^\Pi\subset\sK^\lax(F)$ of $\Pi$-fixed points, and the corestriction $\sK(F)\to\sK^\lax(F)^\Pi$ is an equivalence.
\end{lemma}
\begin{proof}
    Recall that the propagation unit at a lax $F$-algebra $FB\xleftarrow rE\xrightarrow aB$ is the morphism consisting of the horizontal arrows in the diagram
    \begin{align}
        \begin{tikzcd}[ampersand replacement=\&]
            FB \ar[r, "Fi"]
            \& F(B\sqcup_EFB) \\
            E \ar[u, "r"]\ar[r, "r"]\ar[d, "a"']\ar[dr, phantom, "\ulcorner" very near end]
            \& FB \ar[u, "Fi"']\ar[d] \\
            B \ar[r, "i"']
            \& B\sqcup_EFB
        \end{tikzcd}
        \label{eq:pui}
    \end{align}
    By \cref{prop:embed}, $\sK(F)$ can be identified with the full subcategory of $\sK^\lax(F)$ on the lax $F$-algebras $FB\xleftarrow rE\xrightarrow aB$ where the resolution $r$ is an equivalence.
    In particular, this implies that the pushout morphism $i:B\to B\sqcup_EFB$ is an equivalence (since pushouts preserve equivalences), and thus that the top morphism $Fi:FB\to F(B\sqcup_EFB)$ is an equivalence also.
    This shows that the inclusion $\sK(F)\hookrightarrow\sK^\lax(F)$ indeed factors through $\sK^\lax(F)^\Pi$.
    
    Conversely, if a lax $F$-algebra $FB\xleftarrow rE\xrightarrow aB$ lies in $\sK^\lax(F)^\Pi$, then the middle component $r:E\to FB$ in (\ref{eq:pui}) in particular is an equivalence.
    This implies that the lax $F$-algebra lies in the essential image of $\sK(F)$, showing that the fully faithful inclusion $\sK(F)\to\sK^\lax(F)^\Pi$ is essentially surjective, thus completing the proof.
\end{proof}
\begin{remark}\label{rem:LiK}
    By \cite[Corollary 4.4.2.4]{lurie}, finite cocompleteness follows from assuming $\sK$ has pushouts and an initial object, which is necessary to ensure that the propagation endofunctor is well-defined for the entire category $\sK^\lax(F)$.
    This assumption is not strictly necessary: we can instead choose any full subcategory $\sL\subseteq\sK^\lax(F)$ that contains $\sK(F)$ and has enough pushouts to construct a propagation functor $\Pi:\sL\to\sK^\lax(F)$.
    If $\Pi$ corestricts to an endofunctor on $\sL$, then the above lemma can be adapted to show that the inclusion $\sK(F)\hookrightarrow\sL^\Pi$ is an equivalence.
\end{remark}

\begin{theorem}[Ad\'amek's construction on lax algebras]\label{thm:adlax}
    Let $F:\sK\to\sK$ be an endofunctor on an arbitrary $\infty$-category $\sK$, and fix a lax $F$-algebra $FB\xleftarrow rE\xrightarrow aB$.
    Construct the diagrams $D^\theta : [\theta]\to\sK^{\Delta[1]}$ by transfinite induction, where $\theta$ is an ordinal:
    \begin{itemize}
        \item Take $D^0 : [0] \to \sK^{\Delta[1]}$ to be the diagram that picks out the arrow $a:E\to B$.
            Note that the resolution map provides an arrow $r_0 := r:E\to FB$.
        \item Given $D^\theta : [\theta] \to \sK^{\Delta[1]}$, denote by $E^\theta\to B^\theta$ the arrow of $\sK$ picked out by $D^\theta(\theta)$.
            Suppose we have chosen an arrow $r_\theta : E^\theta\to FB^\theta$.
            Then, define $D^{\theta+1}$ to be the extension of $D^\theta$ that sends the morphism $\theta\leq\theta+1$ in $[\theta+1]$ to the pushout square
            $$
            \begin{tikzcd}
                E^\theta \ar[d]\ar[r, "r_\theta"]\ar[dr, phantom, "\ulcorner" very near end]
                & FB^\theta \ar[d] \\
                B^\theta \ar[r, "i_{\theta+1}"']
                & B^{\theta+1}
            \end{tikzcd}
            $$
            In particular, $D^{\theta+1}(\theta+1)$ picks out the arrow $E^{\theta+1}\to B^{\theta+1}$ where $E^{\theta+1} := FB^\theta$.
            Moreover, choose $r_{\theta+1} := Fi_{\theta+1} : E^{\theta+1}\to FB^{\theta+1}$.
        \item For a limit ordinal $\lambda$, and given $D^\theta$ for all $\theta<\lambda$, let $\underline\lambda = \varcolim_{\theta<\lambda}[\theta]$ be the nerve of the poset of all ordinals $0\leq\xi<\lambda$, so that the provided diagrams induce $D^{<\lambda} : \underline\lambda \to \sK^{\Delta[1]}$.

            Then, define $D^\lambda$ to be a colimit cocone for $D^{<\lambda}$.
            If the colimit point $D^\lambda(\lambda)$ is the arrow $E^\lambda\to B^\lambda$, then the choice of $r_\theta$ for every $\theta<\lambda$ induces a canonical map $r_\lambda : E^\lambda\to FB^\lambda$ by the universal property of $E^\lambda$.
    \end{itemize}
    Suppose for some limit ordinal $\lambda$ that the diagram $D^\lambda:[\lambda]\to\sK^{\Delta[1]}$ is well-defined, and let $E^*\to B^*$ be the arrow picked out by $D^\lambda(\lambda)$ in $\sK^{\Delta[1]}$.
    If the canonical map $r_*:E^*\to FB^*$ induced by the $r_\theta$ chosen in the construction is invertible, then the composite $FB^*\xrightarrow{r_*^{-1}}E^*\to B^*$ defines an action that realises $B^*$ as the free $F$-algebra generated by $FB\gets E\to B$.
\end{theorem}
\begin{remark}\label{rem:adlax}
    Ad\'amek's construction for a lax $F$-algebra $FB\xleftarrow rE\xrightarrow aB$ can be described more explicitly if the construction terminates after countably many steps (that is, $\lambda=\omega$).
    In this case, we are assuming that the pushout squares in the diagram
    $$
    \begin{tikzcd}
        E \ar[r, "r"]\ar[d, "a"']\ar[dr, phantom, "\ulcorner" very near end]
        & FB \ar[r, "Fi_1"]\ar[d]\ar[dr, phantom, "\ulcorner" very near end]
        & F(B\sqcup_EFB) \ar[r, "Fi_2"]\ar[d]\ar[dr, phantom, "\ulcorner" very near end]
        & F(B\sqcup_EF(B\sqcup_EFB)) \ar[r, "Fi_3"]\ar[d]
        & \dots \\
        B \ar[r, "i_1"']
        & B\sqcup_EFB \ar[r, "i_2"']
        & B\sqcup_EF(B\sqcup_EFB) \ar[r, "i_3"']
        & B\sqcup_EF(B\sqcup_EF(B\sqcup_EFB)) \ar[r]
        & \dots
    \end{tikzcd}
    $$
    exist, and moreover that this diagram has a colimit $E^*\to B^*$ (in $\sK^{\Delta[1]}$).
    If the canonical map $E^*\to FB^*$ induced by the top row is invertible, then composing an inverse with the colimit arrow defines an action $FB^*\xrightarrow\sim E^*\to B^*$ that realises $B^*$ as the free $F$-algebra generated by $FB\gets E\to B$.
\end{remark}
\begin{proof}[Proof of \cref{thm:adlax}]
    For every $\theta\leq\lambda$, let $E^\theta\to B^\theta$ denote the arrow picked out by $D^\lambda(\theta)$.
    With the chosen arrows $r_\theta : E^\theta\to FB^\theta$, we obtain lax $F$-algebras $A^\theta := \{FB^\theta\gets E^\theta\to B^\theta\}$, where $A^0$ is the original lax $F$-algebra $A^0 = A := \{FB\gets E\to B\}$.

    Let $\sL$ denote the full subcategory of $\sK^\lax(F)$ spanned by $\sK(F)$ and the lax $F$-algebras $A^\theta$ for $\theta\leq\lambda$.
    Assuming that $D^\lambda$ is well-defined ensures that we have the pushouts in $\sK$ necessary to define the unital propagation functor $\Pi:\sL\to\sK^\lax(F)$ as in \cref{def:prop}.
    Moreover, we have by design that $\Pi(A^\theta) = A^{\theta+1}$ for every $\theta<\lambda$, and $\Pi(A^\lambda) \simeq A^\lambda$ since we assume that the map $r_\lambda=r_*$ is invertible.
    Therefore, $\Pi$ corestricts to an endofunctor $\sL\to\sL$.
    
    For $\theta\leq\lambda$, let $D_A^\theta : [\theta]\to\sL$ denote the diagram constructed in \cref{thm:freefix} with the endofunctor $\Pi$ and the lax $F$-algebra $A$.
    We can see by transfinite induction that the diagrams $D_A^\theta$ are indeed well-defined, and moreover that the lax $F$-algebra $D_A^\theta(\theta)$ is precisely $A^\theta$.
    \begin{itemize}
        \item This is immediate if $\theta = 0$.
        \item Given that $D_A^\theta$ is well-defined and $D_A^\theta(\theta) = A^\theta$, it follows that $D_A^{\theta+1}$ exists and maps $(\theta+1)$ to $A^{\theta+1}$ because $A^{\theta+1} = \Pi(A^\theta)$.
        \item Suppose for a limit ordinal $\xi$ that $D_A^\theta$ is well-defined for all $\theta<\xi$, and $D_A^\theta(\theta) = A^\theta$.
            It would follow that $D_A^\xi$ is well-defined and $D_A^\xi(\xi) = A^\xi$ if we can show that the colimit of $D_A^{<\xi} : \underline\xi \to \sL$ is $A^\xi$.

            To see this, recall that the arrow $E^\xi\to B^\xi$ is defined to be the colimit of $D^{<\xi} : \underline\xi\to\sK^{\Delta[1]}$, and the universal property of $E^\xi$ then canonically induces the map $r_\xi$ from the maps $r_\theta$ for $\theta < \xi$.
            This is precisely how the colimit $\varcolim D_A^{<\xi}$ of lax $F$-algebras is constructed; see \cref{prop:laxcol}.
    \end{itemize}
    By assumption, the lax $F$-algebra $A^\lambda = \{FB^*\gets E^*\to B^*\}$ has an invertible resolution $r_*$, so $A^\lambda = D_A^\lambda(\lambda)$ lies in $\sL^\Pi$.
    Therefore, the conclusion follows from \cref{thm:freefix} and \cref{lem:pfp}.
\end{proof}

\begin{corollary}[Ad\'amek's free algebra construction]\label{cor:adamek}
    Let $F:\sK\to\sK$ be an endofunctor on an $\infty$-category $\sK$, and fix any object $K$ of $\sK$.
    Construct the objects $K^\theta$ and morphisms $i_\theta:K^\theta\to K\sqcup FK^\theta$, for $\theta$ an ordinal, by transfinite induction:
    \begin{itemize}
        \item Take $K^0 := K$ and $i_0:K\to K\sqcup FK$ to be the first coprojection.
        \item Given $K^\theta\to K\sqcup FK^\theta$, we have a pushout square
            $$
            \begin{tikzcd}
                FK^\theta \ar[r, "Fi_\theta"]\ar[d]\ar[dr, phantom, "\ulcorner" very near end]
                & F(K\sqcup FK^\theta) \ar[d] \\
                K\sqcup FK^\theta \ar[r]
                & K\sqcup F(K\sqcup FK^\theta)
            \end{tikzcd}
            $$
            where the vertical arrows are given by second coprojections for the respective coproducts.
            Define $K^{\theta+1} := K\sqcup FK^\theta$ and take $i_{\theta+1} : K^{\theta+1}\to K\sqcup FK^{\theta+1}$ to be the bottom row of the above pushout square.
        \item Given a limit ordinal $\lambda$ and $i_\theta:K^\theta\to K\sqcup FK^\theta=K^{\theta+1}$ for every $\theta<\lambda$, define $K^\lambda := \varcolim_{\theta<\lambda}K^\theta$.
            Since $K^{\theta+1} = K\sqcup FK^\theta$, we obtain a canonical map $i_\lambda:K^\lambda\to K\sqcup FK^\lambda$.
    \end{itemize}
    Suppose for some limit ordinal $\lambda$ that $K^\lambda$ is well-defined, and that the induced map $i_\lambda : K^\lambda\to K\sqcup FK^\lambda$ is invertible.
    Then, the canonical map $FK^\lambda\to K^\lambda$ induced by the coprojections $FK^\theta\to K\sqcup FK^\theta = K^{\theta+1}$ for $\theta<\lambda$ defines an action that realises $K^\lambda$ as the free $F$-algebra generated by the object $K$.
\end{corollary}

\begin{proof}
    Assume first that $\sK$ has an initial object $\emptyset$.
    Then, the result follows from applying \cref{thm:adlax} to the free lax $F$-algebra $FK\gets\emptyset\to K$ generated by the object $K$.

    Now, suppose $\sK$ does not have an initial object.
    Extend $F$ to an endofunctor on $\sK^\triangleleft$ by fixing the cone point, then $\sK(F)$ is a full subcategory of $\sK^\triangleleft(F)$.
    Then, $K^\lambda$ is a free $F$-algebra in $\sK^\triangleleft(F)$ generated by $K$ by the previous paragraph.
    Since $K^\lambda$ lives in $\sK(F)$ as well, it restricts to a free $F$-algebra in $\sK(F)$ generated by $K$ also.
\end{proof}
\begin{remark}\label{rem:adamek-countable}
    As in \cref{rem:adlax}, Ad\'amek's construction of free $F$-algebras can be described more succinctly if the construction terminates after countably many steps.
    In this case, we suppose the coproducts in the diagram
    $$
    \begin{tikzcd}
        & FK \ar[r, "Fi_1"]\ar[d]\ar[dr, phantom, "\ulcorner" very near end]
        & F(K\sqcup FK) \ar[r, "Fi_2"]\ar[d]\ar[dr, phantom, "\ulcorner" very near end]
        & F(K\sqcup F(K\sqcup FK)) \ar[r, "Fi_3"]\ar[d]
        & \dots \\
        K \ar[r, "i_1"']
        & K\sqcup FK \ar[r, "i_2"']
        & K\sqcup F(K\sqcup FK) \ar[r, "i_3"']
        & K\sqcup F(K\sqcup F(K\sqcup FK)) \ar[r, "i_4"']
        & \dots
    \end{tikzcd}
    $$
    exist (note that the top row is obtained from the bottom row by applying $F$).
    If the bottom row has a colimit $K^*$ in $\sK$ that is preserved by the functor $K\sqcup F(-):\sK\to\sK$, then the canonical map $FK^*\to K^*$ induced by the vertical arrows in the above diagram realises $K^*$ as the free $F$-algebra generated by $K$.
\end{remark}

\begin{corollary}[Ad\'amek's initial algebra construction]\label{thm:adamek}
    Let $\emptyset$ be an initial object of $\sK$.
    Construct objects $I^\theta$ and maps $i_\xi : I^\xi \to FI^\theta$ for $\xi\leq\theta$ by transfinite induction:
    \begin{itemize}
        \item Define $I^0 := \emptyset$, with $i_0$ the unique map $I^0\to FI^0$.
        \item Given $i_\theta$, define $I^{\theta+1} := FI^\theta$ and $i_{\theta+1} := Fi_\theta$.
        \item For a limit ordinal $\lambda$, given the maps $i_\theta : I^\theta \to FI^\theta = I^{\theta+1}$ for every $\theta<\lambda$, define $I^\lambda := \varcolim_{\theta<\lambda}I^\theta$.
            Since $FI^\theta = I^{\theta+1}$, we then obtain a canonical map $i_\lambda:I^\lambda\to FI^\lambda$.
    \end{itemize}
    Suppose for some limit ordinal $\lambda$ that $I^\lambda$ is well-defined, and that the induced map $i_\lambda : I^\lambda\to FI^\lambda$ is invertible.
    Then, the pair $(I^\lambda, i_\lambda^{-1})$ defines an initial $F$-algebra.
\end{corollary}
\begin{proof}
    Follows from \cref{cor:adamek} by taking $K=\emptyset$.
\end{proof}
\begin{remark}\label{rem:adamek-countable-initial}
    The corollary implies in particular that if
    \[
        \emptyset \xrightarrow{!}F\emptyset \xrightarrow{F(!)} F^2\emptyset \xrightarrow{F^2(!)} F^3\emptyset \to \dots
    \]
    admits a colimit $I$, and the canonical map $f:I\to FI$ is an equivalence, then $(I,f^{-1})$ is an initial $F$-algebra.
\end{remark}

%% file: src/algebra/fixedpoints.tex
%%!TEX root = ../../article.tex

\subsection{Fixed points of an endofunctor}\label{ssec:fix}

Let $F:\sK\to\sK$ be an endofunctor.
We are interested in the fixed points of $F$, with a specified identification between the object $K$ and its image $FK$.
These fixed points with chosen identifications collect to form an $\infty$-category:

\begin{definition}\label{def:fix}
    Let $\bI$ denote the nerve of the walking isomorphism (that is, the groupoid with two objects and a unique morphism between any pair of objects).
    Then, define the $\infty$-category $\Fix(F)$ as the pullback
    $$
    \begin{tikzcd}
        \Fix(F) \ar[r]\ar[d]\ar[dr, phantom, "\lrcorner" very near start]
        & \sK^{\bI} \ar[d, two heads] \\
        \sK \ar[r, "{(F,\id)}"']
        & \sK\times\sK
    \end{tikzcd}
    $$
    Note that the vertical map on the right is a categorical fibration, as it is induced by an inclusion of simplicial sets, meaning that the pullback is a homotopy pullback of $\infty$-categories.
\end{definition}

Recall that we denote by $\sK(F)$ the $\infty$-category of $F$-algebras, and $\sK_\co(F) = (\sK^\op(F))^\op$ the $\infty$-category of $F$-coalgebras.

\begin{proposition}\label{prop:fix}
    Let $F:\sK\to\sK$ be an endofunctor of an arbitrary $\infty$-category $\sK$.
    Then, the $\infty$-category $\Fix(F)$ is equivalent to the full subcategory of $\sK(F)$ spanned by the $F$-algebras with trivial action.
    Dually, $\Fix(F)$ is also equivalent to the full subcategory of $\sK_\co(F)$ spanned by the $F$-coalgebras with trivial coaction.
\end{proposition}
\begin{proof}
    Let $\bI^\natural=\bI^\sharp$ denote the walking isomorphism as a marked simplicial set, marked at every edge.
    This is fibrant in the cartesian model structure on marked simplicial sets.
    Let $\sK_\co(F)^F$ denote the full subcategory of $\sK_\co(F)$ spanned by the $F$-coalgebras with trivial coaction.
    We will establish that the forgetful functor $\Fix(F)\to\sK_\co(F)$ restricts to an equivalence $\Fix(F)\to\sK_\co(F)^F$.
    The statement regarding $F$-algebras is similar.

    Let $\Delta[1]^\sharp$ denote the simplicial set $\Delta[1]$ marked at all edges.
    Then, either inclusion $\Delta[1]^\sharp\to\bI^\natural$ is marked anodyne.
    In particular, one of these inclusions induces a trivial categorical fibration $\sK^{\bI}=\Map^\flat(\bI^\natural,\sK^\natural)\to\Map^\flat(\Delta[1]^\sharp,\sK^\natural)$ fitting in the diagram
    $$
    \begin{tikzcd}
        \Fix(F) \ar[rr]\ar[ddr, bend right]\ar[dr, dashed]
        && \sK^{\bI} \ar[d, two heads, "\sim" sloped] \\
        & \sK_\co(F)^F \ar[r]\ar[d]\ar[dr, phantom, "\lrcorner" very near start]
        & \Map^\flat(\Delta[1]^\sharp,\sK^\natural) \ar[d, two heads] \\
        & \sK \ar[r, "{(\id,F)}"']
        & \sK\times\sK
    \end{tikzcd}
    $$
    The perimeter commutes for a suitably chosen inclusion $\Delta[1]^\sharp\to\bI^\natural$, and is the definitional pullback square for $\Fix(F)$.
    As $\sK^\natural$ is marked at the equivalences, the inner square is also a pullback square.

    All of the vertical maps on the right are categorical fibrations, so the pullback diagrams are both homotopy pullback diagrams in the Joyal model structure on simplicial sets.
    In particular, since the corresponding objects of the two pullback diagrams are connected by categorical equivalences, it follows that the induced map $\Fix(F)\to\sK_\co(F)^F$ is a categorical equivalence as well.
\end{proof}
\begin{corollary}\label{cor:itfp}
    Let $F:\sK\to\sK$ be an endofunctor of an arbitrary $\infty$-category $\sK$.
    Then, any initial $F$-algebra defines an initial object in $\Fix(F)$.
    Dually, any terminal $F$-coalgebra defines a terminal object in $\Fix(F)$.
\end{corollary}
\begin{proof}
    Follows by combining \cref{prop:fix} with \cref{lem:lambek}.
\end{proof}

For a limit ordinal $\lambda$, denote by $\underline\lambda$ the nerve of the poset category of all ordinals $0\leq\xi<\lambda$.
Then, define a $\lambda$-\emph{sequence} in $\sK$ to be a functor $\underline\lambda\to\sK$.
\begin{definition}\label{def:lseq}
    Let $F:\sK\to\sK$ be an endofunctor of an $\infty$-category $\sK$, and fix a limit ordinal $\lambda$.
    Say that the pair $(\sK,F)$ is \emph{compatible with $\lambda$-sequences} if:
    \begin{itemize}
        \item $\sK$ is closed under colimits of $\theta$-sequences for all limit ordinals $0<\theta\leq\lambda$, and
        \item $F$ preserves colimits of $\lambda$-sequences.
    \end{itemize}
\end{definition}
Note that $(\sK,F)$ is automatically a compatible with $\lambda$-sequences for some limit ordinal $\lambda$ if $\sK$ is locally presentable, and $F$ is an accessible functor.
The remainder of this section is dedicated to providing a cursory study of $\Fix(F)$ in the case where $(\sK,F)$ is compatible with $\lambda$-sequences for some limit ordinal $\lambda$.

The key observation is the following:

\begin{theorem}\label{thm:fixloc}
    Suppose $(\sK,F)$ is compatible with $\lambda$-sequences for some limit ordinal $\lambda$.
    Then, the fully faithful inclusion $\Fix(F)\hookrightarrow\sK_\co(F)$ admits a left adjoint, realising $\Fix(F)$ as a reflective localisation of the $\infty$-category of $F$-coalgebras.
\end{theorem}
\begin{proof}
    Note that $F$ induces a unital endofunctor on $\sK_\co(F)$.
    Explicitly, the image under $F$ of an $F$-coalgebra is another $F$-coalgebra, and we have a unit $\eta:\Id\Rightarrow F$ whose component on any coalgebra $(C,\nu)$ is its coaction $\eta_C = \nu$.
    Moreover, the full subcategory $\sK_\co(F)^F$ of the $F$-fixed points in $\sK_\co(F)$ is precisely the full subcategory of $F$-coalgebras with trivial coaction by definition.
    By \cref{prop:fix}, it therefore suffices to show that the inclusion $\sK_\co(F)^F\hookrightarrow\sK_\co(F)$ admits a left adjoint.

    By \cref{prop:cocol}, colimits of $\theta$-sequences in $\sK_\co(F)$ for $\theta\leq\lambda$ exist and are computed on the underlying objects in $\sK$.
    In particular, for any $F$-coalgebra $(C, \nu)$, the colimit $I_C$ of the $\lambda$-sequence
    \[
        C \xrightarrow\nu FC \xrightarrow{F\nu} F^2C \xrightarrow{F^2\nu}F^3C \to \dots
    \]
    in $\sK_\co(F)$ exists, with coaction given by the canonical map $\varcolim_nF(F^nC) \to F(\varcolim_nF^nC)$.
    Since $F$ preserves colimits of $\lambda$-sequences, this coaction is an equivalence.
    By \cref{thm:freefix}, it follows that the functor $\Hom_{\sK_\co(F)}(C, -)| : \sK_\co(F)^F \to \sS$ is corepresentable.
    Since this is true for any $F$-coalgebra $(C, \nu)$, it follows that the inclusion $\sK_\co(F)^F\hookrightarrow\sK_\co(F)$ admits a left adjoint, as desired.
\end{proof}

\begin{definition}\label{def:Flocal}
    Let $(\sK,F)$ be compatible with $\lambda$-sequences for some limit ordinal $\lambda$.
    Denote the left adjoint of the inclusion $\Fix(F)\hookrightarrow\sK_\co(F)$ by $I_{(-)}:\sK_\co(F)\to\Fix(F)$.
    Call an $F$-coalgebra homomorphism $F$-\emph{local} if its image under $I_{(-)}$ is an equivalence in $\Fix(F)$.
\end{definition}
\begin{example}\label{ex:coloc}
    If $(C,\nu)$ is an $F$-coalgebra, then the coaction $\nu$ trivially defines an $F$-coalgebra homomorphism $(C,\nu)\to(FC,F\nu)$.
    As $F$-coalgebra homomorphisms, all coactions are $F$-local.
\end{example}

We now provide a complete characterisation of the $F$-local morphisms:

\begin{proposition}\label{prop:llift}
    Suppose $(\sK,F)$ is compatible with $\lambda$-sequences for some limit ordinal $\lambda$.
    Let $I_{(-)}:\sK_\co(F)\to\Fix(F)$ denote the left adjoint to the inclusion.
    Then, an $F$-coalgebra homomorphism $F:(C,\nu)\to(D,\mu)$ is $F$-local if and only if there exists a morphism $s:D\to I_C$ such that
    $$
    \begin{tikzcd}
        C \ar[r, "\varphi"]\ar[d]
        & D \ar[dl, dashed, "s"]\ar[d] \\
        I_C \ar[r, "I_\varphi"']
        & I_D
    \end{tikzcd}
    $$
    commutes up to homotopy.
\end{proposition}
\begin{proof}
    If $\varphi$ is $F$-local, then $I_C\to I_D$ is an equivalence, which allows us to construct the morphism $s$.
    Conversely, any such morphism $s:D\to I_C$ induces a morphism $s_*:I_D\to I_{I_C}\simeq I_C$ fitting in the diagram
    $$
    \begin{tikzcd}
        I_C \ar[r, "\varphi_*"]\ar[d, "\sim"' sloped]
        & I_D \ar[d, "\sim" sloped]\ar[dl, "s_*"'] \\
        I_C \ar[r, "\varphi_*"']
        & I_D
    \end{tikzcd}
    $$
    showing that $s_*$ is a weak inverse of $\varphi_*$.
\end{proof}

\begin{corollary}\label{cor:llift}
    Suppose $(\sK,F)$ is compatible with $\lambda$-sequences for some limit ordinal $\lambda$.
    Let $\varphi:(C,\nu)\to(D,\mu)$ be an $F$-coalgebra homomorphism.
    If there exists a morphism $s:D\to FC$ in $\sK$ such that the diagram
    $$
    \begin{tikzcd}
        C \ar[r, "\varphi"]\ar[d, "\nu"']
        & D \ar[d, "\mu"]\ar[dl, dashed, "s"] \\
        FC \ar[r, "F\varphi"']
        & FD
    \end{tikzcd}
    $$
    commutes up to homotopy, then $\varphi$ is $F$-local.
\end{corollary}

\begin{remark}\label{rem:coloc}
    By \cref{prop:fix}, we have a completely dual theory as well.
    In particular, if $\sK$ is closed under limits of inverse $\lambda$-sequences (that is, functors $\underline\lambda^\op\to\sK$), and $F$ preserves these limits, then $\Fix(F)$ is a coreflective subcategory of $\sK(F)$; that is, the fully faithful inclusion admits a right adjoint $T_{(-)} : \sK(F)\to\Fix(F)$.

    If we refer to an algebra homomorphism as $F$-\emph{colocal} if its image under $T_{(-)}$ is an equivalence, then we have in particular that an algebra homomorphism $\varphi:(A,\alpha)\to(B,\beta)$ is $F$-colocal whenever we can find a map $s:FB\to A$ such that the diagram
    $$
    \begin{tikzcd}
        FA \ar[r, "F\varphi"]\ar[d, "\alpha"']
        & FB \ar[dl, dashed, "s"]\ar[d, "\beta"] \\
        A \ar[r, "\varphi"']
        & B
    \end{tikzcd}
    $$
    commutes up to homotopy.
\end{remark}

We conclude this section with another universal property satisfied by the free fixed point $I_C$ associated to an $F$-coalgebra $C$ realised by studying $F$-algebras \emph{relative} to $C$.
More specifically, let $(C,\nu)$ be an $F$-coalgebra.
Then, we can define an endofunctor $F_C$ on the undercategory $\sK_{C/}$ by the composite
\[
    F_C : \sK_{C/} \xrightarrow{F} \sK_{FC/} \xrightarrow{\nu^!} \sK_{C/}
\]
where $\nu^!$ acts by precomposition with the coaction $\nu:C\to FC$.

\begin{lemma}\label{lem:slice}
    For an $\infty$-category $\sC$ and any functor $p:S\to\sC$, the forgetful functor $\sC_{p/}\to\sC$ creates colimits indexed by weakly contractible simplicial sets.
\end{lemma}
\begin{proof}
    By \cite[Corollary 2.1.2.2]{lurie}, the forgetful functor $\sC_{p/}\to\sC$ is a left fibration.
    Therefore, by the dual of \cite[Proposition 2.4.2.4]{lurie}, the forgetful functor is a cocartesian fibration where every edge of $\sC_{p/}$ is cocartesian.
    In particular, if $K$ is a weakly contractible simplicial set, then any functor $K^\triangleright\to\sC_{p/}$ is a colimit diagram relative to the forgetful functor $\sC_{p/}\to\sC$ by \cite[Proposition 4.3.1.12]{lurie}.
    By \cite[Proposition 4.3.1.5(2)]{lurie}, this means that such a diagram $K^\triangleright\to\sC_{p/}$ is a colimit diagram if and only if its composite $K^\triangleright\to\sC$ is a colimit diagram, as desired.
\end{proof}
\begin{corollary}\label{cor:slice}
    Given an $\infty$-category $\sC$ and an edge $f:x\to y$ in $\sC$, the functor $\sC_{y/}\xrightarrow{f^!}\sC_{x/}$ creates colimits indexed by weakly contractible simplicial sets.
\end{corollary}
\begin{proof}
    We have a commutative triangle
    $$
    \begin{tikzcd}[column sep=small]
        \sC_{y/} \ar[rr, "f^!"]\ar[dr]
        && \sC_{x/} \ar[dl] \\
        & \sC
    \end{tikzcd}
    $$
    where the downward maps create colimits indexed by weakly contractible simplicial sets by \cref{lem:slice}.
\end{proof}

\begin{proposition}\label{prop:fpup}
    Suppose $(\sK,F)$ is compatible with $\lambda$-sequences for some limit ordinal $\lambda$.
    By \cref{thm:fixloc}, denote by $I_{(-)} : \sK_\co(F)\to\Fix(F)$ the left adjoint to the inclusion.
    For any $F$-coalgebra $(C,\nu)$, the free fixed point $I_C$ and the inverse of the induced equivalence $I_C\xrightarrow\sim FI_C$ define an initial $F_C$-algebra in $\sK_{C/}$; that is, an initial object in $\sK_{C/}(F_C)$.
\end{proposition}
\begin{proof}
    By assumption, $F:\sK\to\sK$ preserves colimits of $\lambda$-sequences.
    Therefore, the induced functor $F_C:\sK_{C/}\to\sK_{C/}$ preserves colimits of $\lambda$-sequences by \cref{cor:slice}.
    In particular, using \cref{thm:adamek}, we can construct the initial $F_C$-algebra by Ad\'amek's construction.
    The initial object in $\sK_{C/}$ is given by the identity on $C$, so Ad\'amek's construction builds the $\lambda$-sequence
    \[
        C \xrightarrow\nu FC \xrightarrow{F\nu} F^2C \xrightarrow{F^2\nu} F^3C \to \dots
    \]
    in $\sK_{C/}$.
    By \cref{lem:slice}, this colimit can be calculated in $\sK$, which is precisely the colimit defining the free fixed point $I_C$ in \cref{thm:fixloc}.
\end{proof}

%% file: src/womcats.tex
%%!TEX root = ../article.tex

\section{Constructing higher categories}\label{sec:womcat}

In \cref{ssec:operad}, we review the necessary machinery to functorially associate to every monoidal $\infty$-category $\sV$ an $\infty$-category $\sV\Cat$ of categories weakly enriched in $\sV$, following \cite{gepner-haugseng}.
In \cref{ssec:endo}, we restrict our attention to symmetric monoidal $\infty$-categories $\sV$, which endows $\sV\Cat$ with canonical symmetric monoidal structure.
We show in \cref{thm:wclim} that enrichment preserves well-behaved limits of diagrams indexed by weakly contractible simplicial sets, laying the necessary groundwork to apply Ad\'amek's construction for $\Cat_{(\infty, \infty)}$.

With the enrichment endofunctor defined, we can inductively produce the $\infty$-category $\Cat_{(n,r)}$ of $(n, r)$-categories for $-2\leq n\leq\infty$ and finite $0\leq r\leq n+2$ as follows:

\begin{definition}\label{def:nrcat}
    Let $\bS^m := \partial\Delta[m+1]$ for $m\geq-1$.
    For $-2\leq n\leq\infty$, define $\Cat_{(n,0)} := \Grpd_n$ to be the full subcategory of $\sS$ spanned by the spaces that are local with respect to the maps $\bS^m\to*$ for $m>n$.
    Note in particular that $\Grpd_\infty = \sS$.

    For $-2\leq n\leq\infty$ and finite $0\leq r\leq n+2$, we may now proceed inductively and define $\Cat_{(n+1, r+1)} := (\Cat_{(n,r)})\Cat$ to be the symmetric monoidal $\infty$-category of categories enriched in the cartesian symmetric monoidal $\infty$-category $\Cat_{(n,r)}$.
\end{definition}
\begin{remark}\label{rem:nrcat}
    \cite[Proposition 6.1.7, Theorem 6.1.8]{gepner-haugseng} imply for $n\leq n'$ and $r\leq r'$ that the inclusion $\Cat_{(n,r)}\hookrightarrow\Cat_{(n',r')}$ exhibits $\Cat_{(n, r)}$ as the localisation with respect to the maps
    \begin{itemize}
        \item $\Sigma^r\bS^j\to\Sigma^r*$ for $n-r<j\leq n'-r$, and
        \item $\Sigma^k[1]\to\Sigma^k*$ for $r\leq k<r'$
    \end{itemize}
    where $[1]$ is the walking $1$-morphism, and $\Sigma X$ is the higher category with two objects $\bot, \top$ with no nontrivial endomorphisms, and $\Hom_{\Sigma X}(\bot, \top) = X$ (see \cite[Definition 4.3.21]{gepner-haugseng}).
\end{remark}

\begin{lemma}\label{lem:adjs}
    For every $-2\leq n\leq\infty$ and all finite $0\leq r\leq r'\leq n+2$, the inclusion $\Cat_{(n,r)}\hookrightarrow\Cat_{(n,r')}$ admits a left adjoint $\pi$ and a right adjoint $\kappa$, both of which preserve products.
\end{lemma}
\begin{proof}
    By \cite[Theorem 4.4.7]{gepner-haugseng}, the $\infty$-category $\Cat_{(\infty,1)}$ as defined above coincides with the $\infty$-category of complete Segal spaces (and thus coincides with the usual $\infty$-category $\Cat_\infty$ of quasicategories).
    In particular, it follows that the inclusion $\Cat_{(\infty,0)}\hookrightarrow\Cat_{(\infty,1)}$ admits both a left adjoint $\pi:\Cat_{(\infty,1)}\to\Cat_{(\infty,0)}$ from the localisation described in \cref{rem:nrcat}, and a right adjoint $\kappa:\Cat_{(\infty,1)}\to\Cat_{(\infty,0)}$ given by sending an $\infty$-category to its maximal sub-$\infty$-groupoid (that is, the underlying space of objects of the complete Segal space).
    Note that both adjoints preserve products.
    Moreover, $\pi$ and $\kappa$ restrict to functors $\Cat_{(n,1)}\to\Cat_{(n,0)}$ for every $-2\leq n\leq\infty$.
    Therefore, the conclusion follows from \cite[Proposition 5.7.17]{gepner-haugseng} by iteratively applying enrichment to the above adjoints $\pi$ and $\kappa$.
\end{proof}

We can now define the following $\infty$-categories of $(\infty,\infty)$-categories:
\begin{definition}\label{def:womcat}
    Define the $\infty$-category of \emph{$(\infty, \infty)$-categories} to be the limit
    \[
        \Cat_{(\infty,\infty)} &:= \varlim\left(\dots\to\Cat_{(\infty,3)}\xrightarrow\kappa\Cat_{(\infty,2)}\xrightarrow\kappa\Cat_{(\infty,1)}\xrightarrow\kappa\Cat_{(\infty,0)}\right)
    \]
    Similarly, define the $\infty$-category of \emph{$\omega$-categories} to be the limit
    \[
        \Cat_\omega &:= \varlim\left(\dots\to\Cat_{(\infty,3)}\xrightarrow\pi\Cat_{(\infty,2)}\xrightarrow\pi\Cat_{(\infty,1)}\xrightarrow\pi\Cat_{(\infty,0)}\right)
    \]
\end{definition}

%Given that enrichment preserves limits of diagrams indexed by weakly contractible simplicial sets, it follows that the canonical functors $(\Cat_\omega)\Cat\to\Cat_\omega$ and $(\Cat_{(\infty,\infty)})\Cat\to\Cat_{(\infty,\infty)}$ are equivalences of $\infty$-categories.
%In \cref{ssec:womcat}, we will show that $\Cat_\omega$ and $\Cat_{(\infty,\infty)}$ are suitably universal with this property, establishing the main result \cref{thm:A} of our paper.

\begin{remark}\label{rem:ciwe}
    Intuitively, an $\omega$-category $\cD$ is a formal sequence $(\dots,\pi_{\leq2}\cD,\pi_{\leq1}\cD,\pi_{\leq0}\cD)$ of higher categories such that $\pi_{\leq n}\cD$ is the $(\infty, n)$-category obtained from $\pi_{\leq n+1}\cD$ by formally inverting the $(n+1)$-morphisms.
    On the other  hand, an $(\infty, \infty)$-category $\cC$ is a formal sequence $(\dots,\kappa_{\leq2}\cC,\kappa_{\leq1}\cC,\kappa_{\leq0}\cC)$ of higher categories such that $\kappa_{\leq n}\cC$ is obtained from $\kappa_{\leq n+1}\cC$ as the maximal sub-$(\infty, n)$-category.
    
    In any higher category, the class of equivalences (of any dimension) must satisfy the following saturation condition:
    \begin{itemize}
        \item[$(*)$] If $f:x\to y$ is a $k$-morphism such that there exists a $k$-morphism $g:y\to x$ and equivalences $fg\simeq\id_y$ and $gf\simeq\id_x$, then $f$ is an equivalence.
    \end{itemize}
    This condition uniquely determines a class of equivalences in any $(n, r)$-category so long as $r$ is finite, but this is not true in general.
    Instead, the saturation condition leads to two natural classes of morphisms.
    \begin{itemize}
        \item Define the class of \emph{inductive equivalences} to be the class generated inductively by the saturation condition $(*)$, starting with the assertion that identity $k$-morphisms are equivalences.
        \item On the other hand, define the class of \emph{reversible morphisms} to be the class generated coinductively by the saturation condition $(*)$.
            Note that all inductive equivalences are reversible.
    \end{itemize}
    Taking the equivalences to be precisely the reversible morphisms recovers the notion of pseudo-invertibility used in \cite{cheng}, and implies that any such higher category wherein every higher morphism admits a dual is, in fact, an $\infty$-groupoid.
    This may be too strong of a notion of equivalence in some cases: for instance, the $(\infty, \infty)$-category $\mathbf{Cob}_\infty$ of cobordisms---determined by its truncations $\kappa_{\leq n}\mathbf{Cob}_\infty\simeq\mathbf{Cob}_n$ being the usual $(\infty, n)$-category of cobordisms for each $n\geq0$---has duals for every higher morphism, but is not an $\infty$-groupoid.
    
    One might expect that the equivalences in an $(\infty, \infty)$-category are precisely the inductive equivalences, whereas the equivalences in an $\omega$-category are precisely the reversible morphisms, but the truth is more subtle.
    In an $\omega$-category, the equivalences are precisely those higher morphisms $f$ which can be reversed up to arbitrarily higher dimension.%; see Remark~\ref{rem:omegawe}.
    In particular, equivalences in an $\omega$-category necessarily include the reversible morphisms.
    However, these two notions do not coincide; Henry and Loubaton construct an explicit example in \cite[Construction 4.29]{henry-loubaton}.
\end{remark}

\begin{remark}\label{rem:ccc}
    It follows by construction that the symmetric monoidal structure on $\Cat_{(\infty, \infty)}$ is cartesian.
    We will see in \cref{prop:inal} that this $\infty$-category is presentably symmetric monoidal, meaning in particular that $\Cat_{(\infty, \infty)}$ is a cartesian closed $\infty$-categories.
    It then follows from \cite[Corollary 7.4.10]{gepner-haugseng} that $\Cat_{(\infty, \infty)}$ is canonically self-enriched.
    Since we will have an equivalence $(\Cat_{(\infty, \infty)})\Cat\simeq\Cat_{(\infty, \infty)}$, this provides a construction of the large $(\infty, \infty)$-category of $(\infty, \infty)$-categories.
\end{remark}

\input{src/womcats/operads}

\input{src/womcats/endofunctor}

\input{src/womcats/omegacat}

\input{src/womcats/noeth}

%% file: src/womcats/operads.tex
%%!TEX root = ../../article.tex

\subsection{Symmetric and nonsymmetric \texorpdfstring{$\infty$}{oo}-operads}\label{ssec:operad}

In this section, we provide a brief overview of the necessary details regarding (generalised) nonsymmetric $\infty$-operads in the sense of \cite{gepner-haugseng} and symmetric $\infty$-operads in the sense of \cite{lurie:algebra}.
The author claims no originality in the content below.

\begin{definition}\label{def:catpat}
    A \emph{categorical pattern} $\frakP$ in the sense of \cite[Definition 3.2.1]{gepner-haugseng} consists of an $\infty$-category $\sC$, a family of diagrams $p_\alpha:K_\alpha^\triangleleft\to\sC$, and a marking on $\sC$ such that every edge in $K_\alpha^\triangleleft$ is sent to a marked edge of $\sC$ via $p_\alpha$.

    A map of categorical patterns from $\frakP=(\sC,\{p_\alpha\}) $to $\frakP'=(\sC',\{p'_\beta\})$ is a map $f:\sC\to\sC'$ of marked simplicial sets such that for every index $\alpha$, there is an index $\beta$ such that $f\circ p_\alpha=p'_\beta$.
\end{definition}

\begin{theorem}\cite[Theorem B.0.20, Proposition B.2.9]{lurie:algebra}\label{thm:catpat}
    For a categorical pattern $\frakP=(\sC,\{p_\alpha:K_\alpha^\triangleleft\to\sC\})$, there is a unique left proper combinatorial simplicial model structure on $\sSet^+_{/\sC}$ such that the cofibrations are the morphisms of $\sSet^+$ whose underlying maps of simplicial sets are monomorphisms, and whose fibrant objects are those $\pi:X\to\sC$ such that
    \begin{itemize}
        \item $\pi$ is an inner fibration,
        \item every marked edge of $\sC$ admits a $\pi$-cocartesian lift in $X$, and these $\pi$-cocartesian lifts of marked edges of $\sC$ are precisely the marked edges of $X$,
        \item for every index $\alpha$, the pullback $\pi_\alpha:X\times_\sC K_\alpha^\triangleleft\to K_\alpha^\triangleleft$ along $p_\alpha$ is the cocartesian fibration associated to a limit cone $K_\alpha^\triangleleft\to\Cat_\infty$
        \item for every index $\alpha$ and any section $s:K_\alpha^\triangleleft\to X\times_\sC K_\alpha^\triangleleft$ of $\pi_\alpha$, the composite $K_\alpha^\triangleleft\xrightarrow sX\times_\sC K_\alpha^\triangleleft\xrightarrow{\pi_\alpha}X$ is a $\pi$-limit cone.
    \end{itemize}
    Denote this model structure by $\sSet^+_\frakP$.

    Moreover, given a map $f:\frakP\to\frakP'$ of categorical patterns, composition with $f$ induces a left Quillen functor $f_!:\sSet^+_\frakP\to\sSet^+_{\frakP'}$.
\end{theorem}
\begin{remark}\label{rem:catpat}
    The model structure induced by the trivial categorical pattern on $\Delta[0]$ (with the unique marking, and no diagrams are chosen) on $\sSet^+$ has as fibrant objects the $\infty$-categories marked at their equivalences.
    Moreover, the $\infty$-category induced by $\sSet^+$ is equivalent to the $\infty$-category $\Cat_\infty$ of $\infty$-categories.

    By \cite[Remark B.2.5]{lurie:algebra}, every model category $\sSet^+_\frakP$ induced by a categorical pattern $\frakP$ is canonically enriched over $\sSet^+$ endowed with the above model structure.
    In particular, if $X$ is $\frakP$-fibrant, then $\Hom_\frakP(-,X):(\sSet^+_\frakP)^\op\to\sSet^+$ is a right Quillen functor.
\end{remark}

\begin{definition}\label{def:nonsym}
    Recall that a morphism $\phi:[n]\to[m]$ is \emph{inert} if it is a subinterval inclusion, meaning $\phi(i) = \phi(0) + i$ for every $0\leq i\leq n$.
    We then construct the following $\infty$-categories:
    \begin{itemize}
        \item Let $\frakO^\gen$ denote the categorical pattern on $\bDelta^\op$ obtained by marking $\bDelta^\op$ at the inert morphisms, and choosing as diagrams the subcategory inclusions $G_{[n]}^\triangleleft\to\bDelta^\op$ for $n\geq0$, where $G_{[n]}^\triangleleft$ is spanned by the objects $[0]$, $[1]$, $[n]$, and the inert morphisms between them.
            Then, define the $\infty$-category $\Opd_\infty^\gen$ of \emph{generalised $\infty$-operads} to be the $\infty$-category associated to the model category $\sSet^+_{\frakO^\gen}$.
        \item Let $\frakO$ denote the categorical pattern on $\bDelta^\op$ obtained by marking the inert morphisms, and choosing as diagrams the subcategory inclusions $K_{[n]}^\triangleleft\to\bDelta^\op$ for $n\geq0$, where $K_{[n]}^\triangleleft$ is spanned by the inert morphisms $[1]\to[n]$.
            Then, define the $\infty$-category $\Opd_\infty$ of \emph{$\infty$-operads} to be the $\infty$-category associated to the model category $\sSet^+_\frakO$.
        \item Let $\frakM$ denote the categorical pattern on $\bDelta^\op$ obtained by marking all morphisms, and taking the diagrams $K_{[n]}^\triangleleft\to\bDelta^\op$ as above.
            Then, define the $\infty$-category $\Mon_\infty$ of \emph{monoidal $\infty$-categories} as the $\infty$-category associated to the model category $\sSet^+_\frakM$.
    \end{itemize}
    The identity functor on $\bDelta^\op$ induces maps of categorical patterns $\frakO^\gen\to\frakO\to\frakM$, and thus left Quillen functors $\sSet^+_{\frakO^\gen}\to\sSet^+_\frakO\to\sSet^+_\frakM$.
    Therefore, we have adjunctions of $\infty$-categories
    $$
    \begin{tikzcd}[column sep=small]
        \Opd_\infty^\gen \ar[rr, bend left=15, "L_\gen"]
        & \bot
        & \Opd_\infty \ar[rr, bend left=15]\ar[ll, bend left=15]
        & \bot
        & \Mon_\infty \ar[ll, bend left=15, "{(-)^\otimes}"]
    \end{tikzcd}
    $$
\end{definition}
\begin{remark}\label{rem:nonsym}
    The morphisms in $\Mon_\infty$ are the strong monoidal functors between monoidal $\infty$-categories.
    Let $\Mon_\infty^\lax$ denote the full subcategory of $\Opd_\infty$ spanned by the image of $(-)^\otimes:\Mon_\infty\to\Opd_\infty$, then $\Mon_\infty^\lax$ is the $\infty$-category of monoidal $\infty$-categories and \emph{lax} monoidal functors.
\end{remark}

\begin{definition}\label{def:sym}
    Let $\bGamma^\op$ denote the category of finite pointed sets, generated by the representatives $\langle n\rangle := \{*,1,\dots,n\}$ for $n\geq0$.
    Call a map $\phi:\langle n\rangle\to\langle m\rangle$ \emph{inert} if every $1\leq j\leq m$ in $\langle m\rangle$ is the image of a unique element of $\langle n\rangle$.
    We now construct the following $\infty$-categories:
    \begin{itemize}
        \item Let $\frakO^\Sigma$ denote the categorical pattern on $\bGamma^\op$ obtained by marking $\bGamma^\op$ at the inert morphisms, and choosing as diagrams the subcategory inclusions $K_{\langle n\rangle}^\triangleleft\to\bGamma^\op$ for $n\geq0$, where $K_{\langle n\rangle}^\triangleleft$ is spanned by the inert maps $\langle1\rangle\to\langle n\rangle$.
            Then, define the $\infty$-category $\Opd_\infty^\Sigma$ of \emph{symmetric $\infty$-operads} to be the $\infty$-category associated to the model category $\sSet^+_{\frakO^\Sigma}$.
        \item Fix a symmetric $\infty$-operad $\pi:\sO\to\bGamma^\op$ (that is, an object of $\Opd_\infty^\Sigma$).
            Let $\frakM_\sO$ denote the categorical pattern on $\sO$ obtained by marking all morphisms, and choosing as diagrams all functors $K_{\langle n\rangle}^\triangleleft\to\sO$ for $n\geq0$ that associate to each inert $\langle1\rangle\to\langle n\rangle$ in $K_{\langle n\rangle}^\triangleleft$ a $\pi$-cocartesian lift in $\sO$.
            Then, define the $\infty$-category $\Mon_\sO$ of \emph{$\sO$-monoidal categories} to be the $\infty$-category associated to the model category $\sSet^+_{\frakM_\sO}$.
    \end{itemize}
\end{definition}
\begin{remark}\label{rem:s2ns}
    We have from \cite[Construction 4.1.2.9]{lurie:algebra} a functor $c:\bDelta^\op\to\bGamma^\op$ defined by sending $[n]$ to the set of partitions (``cuts'') of $[n]$ into at most two contiguous pieces.
    This functor defines a map of categorical patterns $c:\frakO\to\frakO^\Sigma$, and thus an adjunction $c_!:\Opd_\infty\rightleftarrows\Opd_\infty^\Sigma:c^*$, where the right adjoint is the forgetful functor.
\end{remark}
\begin{example}\label{ex:symmon}
    Consider the \emph{commutative $\infty$-operad} $\E_\infty$ given by the identity $\bGamma^\op\to\bGamma^\op$, and let $\SymMon_\infty := \Mon_{\E_\infty}$ denote the $\infty$-category of \emph{symmetric monoidal $\infty$-categories}.
    By definition, $\SymMon_\infty$ is the $\infty$-category associated to the model category given by the categorical pattern obtained by marking $\bGamma^\op$ at all edges, and taking as diagrams the subcategory inclusions $K_{\langle n\rangle}^\triangleleft\to\bGamma^\op$ for $n\geq0$.
    In particular, the identity functor on $\bGamma^\op$ induces a map of categorical patterns $\frakO^\Sigma\to\frakM_{\E_\infty}$ and thus an adjunction $\Opd_\infty^\Sigma\rightleftarrows\SymMon_\infty:(-)^\otimes$, where the right adjoint is a (non-full) inclusion.
\end{example}

\begin{definition}\label{def:opalg}
    For a monoidal $\infty$-category $\sV^\otimes$ (viewed as an $\infty$-operad), we have from \cref{rem:catpat} a right Quillen functor $\Hom_\frakO(-,\sV^\otimes) : (\sSet^+_\frakO)^\op\to\sSet^+$.
    Denote the associated functor of $\infty$-categories by $\Alg_{(-)}(\sV) : \Opd_\infty^\op\to\Cat_\infty$.
    In particular, for any $\infty$-operad $\sO$, we have an $\infty$-category $\Alg_\sO(\sV)$ of $\sO$-algebras in $\sV$.
    Denote the cartesian fibration associated to $\Alg_{(-)}(\sV)$ by $\Alg(\sV)\to\Opd_\infty$.
    This is the \emph{algebra fibration} associated to $\sV$.
\end{definition}
\begin{remark}\label{rem:opalg}
    We can also define a symmetric analogue of operad algebras: for a symmetric monoidal $\infty$-category $\sV^\otimes$ (viewed as a symmetric $\infty$-operad), we also have a right Quillen functor $\Hom_{\frakO^\Sigma}(-,\sV^\otimes) : (\sSet^+_{\frakO^\Sigma})^\op\to\sSet^+$ that induces a functor $\Alg^\Sigma_{(-)}(\sV) : (\Opd_\infty^\Sigma)^\op\to\Cat_\infty$.
\end{remark}
\begin{example}\label{ex:mon=alg}
    Any $\infty$-category $\sC$ with finite products induces a symmetric monoidal structure $\sC^\times$ under the cartesian product by \cite[Corollary 2.4.1.9]{lurie:algebra}.
    In particular, we have a cartesian symmetric monoidal $\infty$-category $\Cat_\infty^\times$.
    Now, \cite[Remark 2.4.2.6]{lurie:algebra} establishes an equivalence $\Alg_\sO^\Sigma(\Cat_\infty^\times)\simeq\Mon_\sO$ for any symmetric $\infty$-operad $\sO$.
    In particular, $\SymMon_\infty \simeq \Alg^\Sigma_{\E_\infty}(\Cat_\infty^\times)$ establishes an equivalence between symmetric monoidal $\infty$-categories and commutative algebra objects in $\Cat_\infty^\times$.
\end{example}
\begin{remark}\label{rem:E1}
    Let $\E_1$ denote the \emph{associative $\infty$-operad} given in \cite[Definition 4.1.1.3]{lurie:algebra}.
    We have by \cite[Corollary 4.3.12]{gepner-haugseng} that $\Mon_\infty \simeq \Mon_{\E_1} \simeq \Alg_{\E_1}^\Sigma(\Cat_\infty^\times)$.
\end{remark}

Let $\Pres_\infty$ denote the subcategory of $\Cat_\infty$ spanned by the locally presentable $\infty$-categories and cocontinuous functors (that is, left adjoints).
Recall from \cite[Remark 4.8.1.6]{lurie:algebra} that $\Pres_\infty$ admits a symmetric monoidal tensor product with the universal property tht cocontinuous functors $\sA\otimes\sB\to\sC$ correspond to ordinary functors $\sA\times\sB\to\sC$ that are cocontinuous in each variable.
Using this monoidal structure, we can define locally presentable analogues of monoidal $\infty$-categories:

\begin{definition}\label{def:presmon}
    For a symmetric $\infty$-operad $\sO$, define $\Mon_\sO^\Pres := \Alg_\sO(\Pres_\infty^\otimes)$ to be the $\infty$-category of \emph{presentably $\sO$-monoidal $\infty$-categories}.
    In particular, define the $\infty$-category of \emph{presentably monoidal $\infty$-categories} to be $\Mon_\infty^\Pres := \Mon_{\E_1}^\Pres$, and the $\infty$-category of \emph{presentably symmetric monoidal $\infty$-categories} to be $\SymMon_\infty^\Pres := \Mon_{\E_\infty}^\Pres$.
\end{definition}

\begin{lemma}\label{lem:monlim}
    All of the forgetful functors in the diagram
    $$
    \begin{tikzcd}
        \SymMon_\infty^\Pres \ar[r]\ar[d]
        & \Mon_\infty^\Pres \ar[r]\ar[d]
        & \Pres_\infty \ar[d] \\
        \SymMon_\infty \ar[r]
        & \Mon_\infty \ar[r]
        & \Cat_\infty
    \end{tikzcd}
    $$
    create limits.
\end{lemma}
\begin{proof}
    Note that \cite[Proposition 3.2.2.1]{lurie:algebra} implies for any symmetric monoidal $\infty$-category $\sV^\otimes$ and any symmetric $\infty$-operad $\sO$ that the forgetful functor $\Alg_\sO^\Sigma(\sV)\to\sV$ creates limits.
    Therefore, all of the functors on the bottom row of the above diagram create limits by noting that $\Mon_\infty \simeq \Alg^\Sigma_{\E_1}(\Cat_\infty^\times)$ and
    \[
        \SymMon_\infty &\simeq \Alg^\Sigma_{\E_\infty}(\Cat_\infty^\times) \simeq \Alg^\Sigma_{\E_\infty\otimes\E_1}(\Cat_\infty^\times) \\
            &\simeq \Alg^\Sigma_{\E_\infty}(\Alg^\Sigma_{\E_1}(\Cat_\infty^\times)^\times) \simeq \Alg^\Sigma_{\E_\infty}(\Mon_\infty^\times)
    \]
    using the closed symmetric monoidal structure on symmetric $\infty$-operads described in \cite[\S2.2.5, \S3.2.4]{lurie:algebra}, and the Dunn Additivity Theorem $\E_\infty\otimes\E_1\simeq\E_1$.
    The fact that the functors on the top row create limits follows similarly.

    That the vertical functors create limits follows by \cite[Proposition 5.5.3.13]{lurie}.
\end{proof}

%% file: src/womcats/endofunctor.tex
%%!TEX root = ../../article.tex

\subsection{The enrichment functor}\label{ssec:endo}

We begin by reviewing the construction in \cite{gepner-haugseng} of the $\infty$-category $\sV\Cat$ of categories enriched in a general monoidal $\infty$-category $\sV^\otimes$.
We then prove in \cref{thm:wclim} that the endofunctor induced by enrichment on (presentably) symmetric monoidal $\infty$-categories preserves a broad class of limits: namely, limits of diagrams indexed by a weakly contractible simplicial set.

\begin{definition}\label{def:algcat}
    For a space $S$, let $\bDelta_S^\op\to\bDelta^\op$ be the cocartesian fibration associated to the functor $\bDelta^\op\to\Cat_\infty$ mapping $[n]\mapsto S^{\times(n+1)}$, where the degeneracies are given by diagonal functors, and faces by projections.
    As described in \cite[\S4.1]{gepner-haugseng}, this construction defines a functor $\bDelta^\op_{(-)}:\sS\to\Opd_\infty^\gen$.

    For a monoidal $\infty$-category $\sV^\otimes$, define the $\infty$-category $\Alg_\cat(\sV)$ of \emph{$\sV$-categorical algebras} as the pullback
    $$
    \begin{tikzcd}
        \Alg_\cat(\sV) \ar[r]\ar[d]\ar[dr, phantom, "\lrcorner" very near start]
        & \Alg(\sV) \ar[d] \\
        \sS \ar[r, "L_\gen\bDelta^\op_{(-)}"']
        & \Opd_\infty
    \end{tikzcd}
    $$
    where the vertical map on the right is the algebra fibration of \cref{def:opalg}.
    In particular, a $\sV$-categorical algebra with space of objects $S$ is precisely a map of generalised $\infty$-operads $\cC:\bDelta_S^\op\to\sV^\otimes$.
\end{definition}

\begin{definition}\label{def:VCat}
    Fix a monoidal $\infty$-category $\sV$.
    For an $\infty$-groupoid $S$, define the \emph{trivial $\sV$-category} $E_S^\sV$ on $S$ to be the composite $\bDelta_S^\op\to\bDelta^\op\xrightarrow{\bfB\1}\sV^\otimes$, where $\bfB\1$ is the delooping of the tensor unit of $\sV^\otimes$ viewed as a monoid (see \cite[Proposition 3.1.18]{gepner-haugseng}).
    In particular, let $E^1 := E_{\{0,1\}}^\sV$ be the walking $\sV$-enriched equivalence, and $E^0 := E_{\{0\}}^\sV = \bfB\1$.

    Say that a $\sV$-categorical algebra $\cC$ is a \emph{$\sV$-category} if it is local with respect to the canonical map $E^1\to E^0$.
    Then, define the $\infty$-category $\sV\Cat$ of $\sV$-categories to be the full subcategory of $\Alg_\cat(\sV)$ spanned by the $\sV$-categories.
    By \cite[Corollary 5.7.6]{gepner-haugseng}, this construction defines a functor $(-)\Cat : \Mon_\infty^\lax\to\Cat_\infty$.
\end{definition}
\begin{remark}\label{rem:ffes}
    By \cite[Theorem 5.6.6]{gepner-haugseng}, the $\infty$-category $\sV\Cat$ is precisely the localisation of $\Alg_\cat(\sV)$ at the $\sV$-functors that are fully faithful and essentially surjective.
\end{remark}

\begin{proposition}\label{prop:sym-mon-enr}\cite[Corollary 5.7.12, Proposition 5.7.16]{gepner-haugseng}
    If $\sV$ is a symmetric monoidal $\infty$-category, then $\sV\Cat$ admits a symmetric tensor product as well.
    If $\sV$ is moreover presentably symmetric monoidal, then so is $\sV\Cat$.
    In particular, the enrichment functor canonically restricts to an endofunctor on $\SymMon_\infty$, and restricts further to an endofunctor on $\SymMon_\infty^\Pres$ also.
\end{proposition}
\begin{remark}
    If $(\cV,\otimes)$ is an ordinary symmetric monoidal category, then the tensor product of $\cV$-enriched categories $\cC$ and $\cD$ is given by taking $\Ob(\cC\otimes\cD) := \Ob(\cC)\times\Ob(\cD)$ for the class of objects, and taking $\Hom_{\cC\otimes\cD}((c,d), (c',d')) := \Hom_\cC(c,c')\otimes\Hom_\cD(d,d')$.
    The tensor product of \cref{prop:sym-mon-enr} is a homotopy-coherent generalisation of this.
\end{remark}

The universal properties asserted in \cref{thm:A} will follow from proving that $\Cat_\omega$ underlies a terminal coalgebra for the enrichment endofunctor on $\SymMon_\infty$, and that $\Cat_{(\infty, \infty)}$ underlies an initial algebra for the enrichment endofunctor on $\SymMon_\infty^\Pres$.
In order to establish these universal properties, it is therefore pertinent that we understand the behaviour of limits under the enrichment endofunctor.

The remainder of this section is dedicated to proving the following result:
\begin{theorem}\label{thm:wclim}
    Suppose $(\sV_i)_{i\in K}$ is a diagram in $\Mon_\infty$ indexed by a weakly contractible simplicial set $K$.
    If each of the functors $F_{i, j} : \sV_i \to \sV_j$ induces a natural equivalence
    $$
    \Hom_{\sV_i}(\1_{\sV_i}, -) \Rightarrow \Hom_{\sV_j}(\1_{\sV_j}, F(-))
    $$
    for all edges $i\to j$ in $K$, then the induced map
    $$
    \left(\varlim_{i\in K}\sV_i\right)\Cat \to \varlim_{i\in K}(\sV_i\Cat)
    $$
    is an equivalence of categories.
\end{theorem}

In order to prove \cref{thm:wclim}, we will piece through the construction of the enrichment functor, and study the limits preserved at each step.

\begin{lemma}\label{lem:pint}
    Let $\sC,\sD$ be $\infty$-categories, and $F:\sC^\op\times\sD\to\Cat_\infty$ a functor.
    For a simplicial set $K$, suppose $\sD$ has all $K$-indexed limits, and that $F_c:\sD\to\Cat_\infty$ preserves these limits for all $c\in\sC_0$.
    Then, the corresponding functor $\sD\to\Cat_{\infty,/\sC}$ preserves $K$-indexed limits.
\end{lemma}
\begin{proof}
    Consider the adjunct functor $\sD\to\Fun(\sC^\op,\Cat_\infty)$.
    Since limits in functor $\infty$-categories are computed pointwise by \cite[Corollary 5.1.2.3]{lurie}, the assumptions of the lemma imply that this adjunct functor preserves $K$-indexed limits.
    The desired functor is the composite $\sD\to\Fun(\sC^\op,\Cat_\infty)\to\Cat_{\infty,/\sC}$ of this adjunct with unstraightening, the latter of which is a right adjoint by \cite[Theorem 3.1.5.1(A0)]{lurie}. 
\end{proof}
\begin{corollary}\label{cor:acspl}
    The functor $\Alg_\cat : \Mon_\infty\to\Cat_{\infty,/\sS}$ preserves all limits.
\end{corollary}
\begin{proof}
    Applying \cref{lem:pint} to $\Hom_{\Opd_\infty}(-,-):\Opd_\infty^\op\times\Opd_\infty\to\Cat_\infty$, we see that the corresponding functor $\Opd_\infty\to\Cat_{\infty,/\Opd_\infty}$ preserves all limits.
    Note that the restriction of this functor to $\Mon_\infty^\lax$ is precisely $\Alg(-):\Mon_\infty^\lax\to\Cat_{\infty,/\Opd_\infty}$.
    Since the inclusion $\Mon_\infty\to\Opd_\infty$ is a right adjoint, it follows that $\Alg : \Mon_\infty\to\Opd_\infty\to\Cat_{\infty,/\Opd_\infty}$ is continuous.
    Observing that $\Alg_\cat$ is recovered as the composite
    \[
        \Mon_\infty\xrightarrow{\Alg}\Cat_{\infty,/\Opd_\infty}\xrightarrow{(L_\gen\bDelta^\op_{(-)})^*}\Cat_{\infty,/\sS}
    \]
    and base change is a right adjoint, the result follows.
\end{proof}

\begin{corollary}\label{cor:acwcl}
    The functor $\Alg_\cat : \Mon_\infty \to \Cat_\infty$ preserves limits of diagrams indexed by weakly contractible simplicial sets.
\end{corollary}
\begin{proof}
    Follows by combining \cref{cor:acspl} with the dual of \cref{lem:slice}.
\end{proof}

\begin{proof}[Proof of \cref{thm:wclim}]
    Let $(\sV_i)_{i\in K}$ be a weakly contractible diagram of monoidal categories where each $F_{i,j} : \sV_i\to\sV_j$ induces a natural equivalence
    $$
    \Hom_{\sV_i}(\1_{\sV_i},-)\Rightarrow\Hom_{\sV_j}(\1_{\sV_j},F(-))
    $$
    of spaces.

    Let $\sV := \varlim_i\sV_i$ in $\Mon_\infty$ with strongly monoidal projections $F_i : \sV \to \sV_i$.
    Then, $\Hom_\sV(\1_\sV, -) \simeq \varlim_i\Hom_{\sV_i}(\1_{\sV_i}, F_i(-))$ is a limit of an essentially constant diagram, by assumption, and therefore each $F_i:\sV\to\sV_i$ induces a natural equivalence $\Hom_\sV(\1_\sV, -) \Rightarrow \Hom_{\sV_i}(\1_{\sV_i}, F_i(-))$.
    
    Recall that $\sV\Cat$ is the full subcategory of $\Alg_\cat(\sV)$ spanned by the objects that are local with respect to the morphism $s_\sV:E^1\to E^0$.
    Note for any space $S$ that the image of the trivial $\sV$-category $E_S^\sV$ under the projection $\sV\to\sV_i$ is precisely $E_S^{\sV_i}$.
    Indeed, since $\sV\to\sV_i$ is strongly monoidal, it preserves the tensor unit of $\sV$, and so the composite $\bDelta^\op_S\to\bDelta^\op\xrightarrow{\bfB\1_\sV}\sV^\otimes\to\sV_i^\otimes$ is equivalent to $E_S^{\sV_i}$.

    Moreover, since $F_i$ induces a natural equivalence $\Hom_{\sV_i}(\1_{\sV_i},-)\Rightarrow\Hom_{\sV_j}(\1_{\sV_j},F(-))$, it follows that the induced functor $F_{i,*} : \Alg_\cat(\sV) \to \Alg_\cat(\sV_i)$ induces natural equivalences $\Hom_\cat(E_S^\sV, -)\Rightarrow\Hom_\cat(E_S^{\sV_i}, F_{i,*}(-))$ for every $S$.
    In particular, each $F_{i,*}$ sends $\sV$-categories to $\sV_i$-categories.
    
    If $\cC$ is a categorical $\sV$-algebra, denote by $\cC_i := F_i(\cC)$ the categorical $\sV_i$-algebra induced by the canonical projection $F_i:\sV\to\sV_i$.
    We have an equivalence $\Alg_\cat(\sV) \simeq \varlim_i\Alg_\cat(\sV_i)$ by \cref{cor:acwcl}.
    Therefore, we have for all $\sC,\cD$ in $\Alg_\cat(\sV)$ that $\Hom_\cat(\cD,\cC) \simeq \varlim_i\Hom_\cat(\cD_i,\cC_i)$.
    In particular, if $\cC$ is a categorical $\sV$-algebra such that $\cC_i$ is a $\sV_i$-category for every $i$, then the homotopy equivalences $\Hom_\cat(E^0_i,\cC_i)\to\Hom_\cat(E^1_i,\cC_i)$ induce a homotopy equivalence $\varlim_i\Hom_\cat(E^0_i,\cC_i)\to\varlim_i\Hom_\cat(E^1_i,\cC_i)$.
    From the above discussion, this is precisely the map $\Hom_\cat(E^0,\cC)\to\Hom_\cat(E^1,\cC)$ induced by $s_\sV$.
    Therefore, $\cC$ is a $\sV$-category.
    
    Altogether, this proves that $\sC\in\Alg_\cat(\sV)$ is a $\sV$-category if and only if every projection $\sC_i\in\Alg_\cat(\sV_i)$ is a $\sV_i$-category.
    In other words, the equivalence of categories $\Alg_\cat(\sV) \simeq \varlim_i\Alg_\cat(\sV_i)$ restricts to an equivalence $\sV\Cat \simeq \varlim_i\left(\sV_i\Cat\right)$.
\end{proof}

%% file: src/womcats/omegacat.tex
%%!TEX root = ../../article.tex

\subsection{\texorpdfstring{$(\infty,\infty)$}{(oo, oo)}-categories}\label{ssec:womcat}

Let $\Alg_\Enr := \SymMon_\infty((-)\Cat)$ denote the $\infty$-category whose objects are pairs $(\sV,\tau)$, where $\sV$ is a symmetric monoidal $\infty$-category and $\tau:\sV\Cat\to\sV$ is a symmetric monoidal functor, and whose morphisms $(\sV,\tau)\to(\sV',\tau')$ are symmetric monoidal functors $\sV\to\sV'$ that respect $\tau$ and $\tau'$ up to homotopy.
Then, define $\Fix_\Enr$ to be the full subcategory of $\Alg_\Enr$ spanned by those pairs $(\sV,\tau)$ where $\tau$ is an equivalence.
Entirely analogously, define $\Alg_\Enr^\Pres := \SymMon_\infty^\Pres((-)\Cat)$ and $\Fix_\Enr^\Pres$ by considering only cocontinuous symmetric monoidal functors between presentably symmetric monoidal $\infty$-categories.

\begin{definition}\label{def:tau}
    For $0\leq n\leq\infty$ and finite $0\leq r\leq n$, define $\tau:\Cat_{(n+1,r+1)}\to\Cat_{(n, r)}$ to be the composite $\Cat_{(n+1,r+1)}\xrightarrow\kappa\Cat_{(n+1,r)}\xrightarrow\pi\Cat_{(n, r)}$, where $\kappa$ is the right adjoint to the inclusion provided by \cref{lem:adjs}, and $\pi$ is a localisation.

    Intuitively, $\tau$ sends an $(n+1,r+1)$-category $\sC$ to the $(n, r)$-category whose $n$-morphisms are equivalence classes of $n$-morphisms in $\sC$.
\end{definition}

\begin{lemma}\label{lem:tau}
    For $0\leq n\leq\infty$ and finite $0\leq r < n$, we have a commutative square
    $$
    \begin{tikzcd}
        \Cat_{(n+1,r+1)} \ar[r, "\kappa"]\ar[d, "\pi"']
        & \Cat_{(n+1, r)} \ar[d, "\pi"] \\
        \Cat_{(n, r+1)} \ar[r, "\kappa"']
        & \Cat_{(n, r)}
    \end{tikzcd}
    $$
\end{lemma}
\begin{proof}
    Since $\kappa : \Cat_{(n+1,r+1)}\to\Cat_{(n+1,r)}$ is the image of $\kappa : \Cat_{(n, r)}\to\Cat_{(n, r-1)}$ under $(-)\Cat$ when $r > 0$, and likewise $\pi : \Cat_{(n+1, r)}\to\Cat_{(n, r)}$ is the image of $\pi : \Cat_{(n, r-1)}\to\Cat_{(n-1, r-1)}$ under $(-)\Cat$ when $n, r > 0$, it suffices to prove that the square
    $$
    \begin{tikzcd}
        \Cat_{(n+1,1)} \ar[r, "\kappa"]\ar[d, "\pi"']
        & \Grpd_{n+1} \ar[d, "\pi"] \\
        \Cat_{(n,1)} \ar[r, "\kappa"']
        & \Grpd_n
    \end{tikzcd}
    $$
    commutes for $n\geq1$, as the cases where $r>0$ then follow by iteratively applying $(-)\Cat$ to the above diagram.
    Note that the localisations are trivial when $n=\infty$, so assume $n<\infty$.

    The localisation functors $\pi$ in this case can be described using the functors $h_n : \Cat_\infty\to\Cat_{(n,1)}$ defined in \cite[Proposition 2.3.4.12]{lurie}.
    Explicitly, for a simplicial set $K$, let $[K,\sC]_n$ be the subset of $\Map(\sk^nK,\sC)$ consisting of restrictions of maps $\sk^{n+1}K\to\sC$.
    Then, the $k$-cells of $h_n\sC$ are homotopy classes of maps in $[\Delta[k],\sC]_n$ relative to $\sk^{n-1}(\Delta[k])$.

    For an $(n+1,1)$-category $\sC$, then we have $\pi\sC = h_n\sC$.
    Note that the $k$-cells of $\kappa h_n\sC$ are given by the $k$-cells of $h_n\sC$ whose edges are all invertible in $h_n\sC$.
    On the other hand, the $k$-cells of $h_n\kappa\sC$ are given by homotopy equivalence classes of maps in $[\Delta[k],\kappa\sC]_n$ relative to $\sk^{n-1}(\Delta[k])$, where $k$-cells of $\kappa\sC$ are the $k$-cells of $\sC$ whose edges are all invertible in $\sC$.
    If $\sC$ is an $(n+1, 1)$-category, then all higher morphisms are invertible, which implies that an edge of $\sC$ is invertible if and only if its image in $h_n\sC$ is invertible.
    Therefore, both $\kappa h_n\sC$ and $h_n\kappa\sC$ describe the same simplicial set: the $k$-cells are homotopy classes of maps in $[\Delta[k],\sC]_n$ whose edges are all invertible in $\sC$.
    This proves that $\kappa\pi\sC=\kappa h_n\sC=h_n\kappa\sC=\pi\kappa\sC$, as desired.
\end{proof}

\begin{proposition}\label{prop:nrwoc}
    For all $0\leq n\leq\infty$ and finite $0\leq r\leq n$, the limit
    $$
    \Enr^\infty(\Cat_{(n, r)}, \tau) := \varlim\left(\cdots\Cat_{(n+2, r+2)}\xrightarrow\tau\Cat_{(n+1, r+1)}\xrightarrow\tau\Cat_{(n, r)}\right)
    $$
    is canonically equivalent to $\Cat_{(\infty, \infty)}$.
\end{proposition}
\begin{proof}
    The truncation functors are colocal by the dual of \cref{ex:coloc}, so we may assume without loss of generality that $r=0$ by induction.
    Observe for $0\leq n<\infty$ that $\Enr^\infty(\Grpd_n,\tau)\simeq\Enr^\infty(\Grpd_{n+1},\tau)$.
    Indeed, \cref{lem:tau} ensures that we have commutativity of the diagram
    $$
    \begin{tikzcd}
        \Cat_{(n+2,1)} \ar[d, "\kappa"']\ar[r, "\pi"] \ar[dd, bend right=60, "\tau"']
        & \Cat_{(n+1,1)} \ar[ddl, "\kappa"]\ar[d, "\kappa"]\ar[dd, bend left=60, "\tau"] \\
        \Grpd_{n+2} \ar[d, "\pi"']
        & \Grpd_{n+1} \ar[d, "\pi"] \\
        \Grpd_{n+1} \ar[r, "\pi"']
        & \Grpd_n
    \end{tikzcd}
    $$
    so the equivalence follows from \cref{rem:coloc}.
    Now, the conclusion follows from the fact that $\Grpd_\infty\simeq\varlim_n\Grpd_n$, noting that the truncation map $\tau:\Cat_{(n,1)}\to\Grpd_n$ reduces to $\kappa$ when $n=\infty$.
\end{proof}

We now prove that $\Cat_{(\infty,\infty)}$ enjoys a universal property dual to that of $\Cat_\omega$ in the presentable setting.

\begin{proposition}\label{prop:inal}
    The $\infty$-category $\Cat_{(\infty,\infty)}$ defines an initial object in $\Alg_\Enr^\Pres$, the $\infty$-category of algebras for the endofunctor $(-)\Cat$ over $\SymMon_\infty^\Pres$.
\end{proposition}
\begin{proof}
    By definition, $\Cat_{(\infty,\infty)} := \Enr^\infty(\Grpd_\infty,\kappa)$ is given as the limit of $\infty$-categories
    \[
        \Cat_{(\infty,\infty)} &:= \varlim\left(\dots\to\Cat_{(\infty,3)}\xrightarrow\kappa\Cat_{(\infty,2)}\xrightarrow\kappa\Cat_{(\infty,1)}\xrightarrow\kappa\Cat_{(\infty,0)}\right)
    \]
    where every $\kappa$ is a right adjoint.
    If $\Pres_\infty^R$ denotes the $\infty$-category of locally presentable $\infty$-categories and right adjoint functors, then \cite[Theorem 5.5.3.18]{lurie} states that the forgetful functor $\Pres_\infty^R\hookrightarrow\Cat_\infty$ creates limits.
    Since \cite[Corollary 5.5.3.4]{lurie} establishes an equivalence $\Pres_\infty \simeq (\Pres_\infty^R)^\op$ by sending cocontinuous functors to their right adjoints, it follows that we can equivalently calculate $\Cat_{(\infty,\infty)}$ as the colimit
    \begin{align}
        \Cat_{(\infty,\infty)} &\simeq \varcolim\left(\Cat_{(\infty,0)}\hookrightarrow \Cat_{(\infty,1)} \hookrightarrow \Cat_{(\infty,2)} \hookrightarrow \dots \right) \label{eq:womcl}
    \end{align}
    in $\Pres_\infty$.

    Now, the forgetful functor $\SymMon_\infty^\Pres := \Alg_{\E_\infty}^\Sigma(\Pres_\infty^\otimes) \to \Pres_\infty$ creates sifted colimits by \cite[Corollary 3.2.3.2]{lurie:algebra}.
    Therefore, we obtain the presentably symmetric monoidal $\infty$-category $\Cat_{(\infty,\infty)}$ as the colimit (\ref{eq:womcl}) computed in $\SymMon_\infty^\Pres$, and this colimit is moreover preserved by $(-)\Cat$ because the corresponding limit of right adjoints is.

    By \cite[Remark 3.1.25]{gepner-haugseng}, the initial object of $\SymMon_\infty^\Pres$ is $\Grpd_\infty$ with its cartesian monoidal structure.
    Therefore, the colimit (\ref{eq:womcl}) in $\SymMon_\infty^\Pres$ is precisely Ad\'amek's construction of an initial algebra for $(-)\Cat$ described in \cref{thm:adamek}.
\end{proof}

We have now proven all of the necessary prerequisites for \cref{thm:A}:

\begin{theorem}\label{thm:main}
    The $\infty$-category $\Cat_{(\infty,\infty)}$ defines an initial object in $\Fix_\Enr^\Pres$.
\end{theorem}
\begin{proof}
    Follows by combining \cref{prop:inal} with \cref{cor:itfp}.
\end{proof}

%% file: src/womcats/noeth.tex
%%!TEX root = ../../article.tex

\subsection{Noetherian \texorpdfstring{$(\infty, \infty)$}{(oo, oo)}-categories}

The crux of the proof of Proposition~\ref{prop:inal} is the observation that $\Cat_{(\infty, \infty)}$ can be realised as an instance of Ad\'amek's initial algebra construction over $\SymMon_\infty^\Pres$.
Nonetheless, the enrichment endofunctor remains canonically-defined in the larger category $\SymMon_\infty$, and it is natural to wonder what the initial algebra of enrichment is in this non-presentable setting.

\begin{proposition}\label{prop:noeth-motivation}
    Ad\'amek's construction of an initial algebra for $(-)\Cat:\SymMon_\infty\to\SymMon_\infty$ does not terminate after $\omega$ steps.
\end{proposition}
\begin{proof}
    The initial object in $\SymMon_\infty$ is the one-object category with its unique symmetric tensor product, which is equivalent to $\Cat_{(-2,0)}$.
    In particular, after $\omega$ steps, Ad\'amek's construction produces the colimit
    \[
        \Cat_{<\omega} := \varcolim\left(\Cat_{(-2,0)} \subseteq \Cat_{(-1, 1)} \subseteq \Cat_{(0, 2)} \subseteq \Cat_{(1, 3)} \subseteq\dots\right) = \bigcup_{\substack{0\leq n<\infty \\ r\geq0}}\Cat_{(n, r)}
    \]
    which is the $\infty$-category of finite-dimensional higher categories.

    However, the objects of $(\Cat_{<\omega})\Cat$ are the $(\infty, \infty)$-categories $\cC$ that are locally finite-dimensional, in the sense that $\Hom_\cC(x, y)$ is a finite-dimensional higher category for all pairs of objects $x, y$ in $\cC$.
    This is a strictly larger $\infty$-category, meaning that $\Cat_{<\omega}\to(\Cat_{<\omega})\Cat$ is not an equivalence.
\end{proof}

To better understand Ad\'amek's construction in this setting, we introduce the following measure of finiteness to $(\infty, \infty)$-categories:

\begin{definition}\label{def:cat-rank}
    We define the \emph{rank} of an $(\infty, \infty)$-category $\sC$ by transfinite induction.
    \begin{itemize}
        \item Say that $\rank\sC < 0$ if and only if $\sC\simeq*$.
        \item For an ordinal $\theta$, say that $\rank\sC < \theta + 1$ if $\rank\Hom_\sC(x, y) < \theta$ for all $x, y\in\sC$.
        \item For a limit ordinal $\lambda$, say that $\rank\sC < \lambda$ if $\rank\sC < \theta$ for some $\theta < \lambda$.
    \end{itemize}
    Say $\rank\sC = \theta$ if $\rank\sC < \theta+1$ but $\rank\sC\not<\theta$.
    Note that the rank of $\sC$ is invariant under equivalence.

    For an ordinal $\theta$, let $\Cat_{<\theta}$ denote the full subcategory of $\Cat_{(\infty, \infty)}$ spanned by those $\sC$ with $\rank\sC < \theta$.
\end{definition}
\begin{remark}\label{rem:cat-rank}
    By \cref{lem:cat-rank} below, if $\rank\sC<\theta$ and $\theta<\theta'$, then also $\rank\sC<\theta'$.
\end{remark}
\begin{example}\label{ex:cat-rank}
    As in \cref{prop:noeth-motivation}, the category $\Cat_{<\omega}$ consists of the finite-dimensional higher categories, and $\Cat_{<\omega+1}$ consists of the locally finite-dimensional higher categories.
\end{example}
\begin{lemma}\label{lem:cat-rank}
    The categories $\Cat_{<\theta}$ can be constructed through transfinite induction:
    \begin{itemize}
        \item $\Cat_{<0} \simeq \Grpd_{-2} \simeq \{*\}$,
        \item $\Cat_{<\theta + 1} \simeq (\Cat_{<\theta})\Cat$; in particular, $\Cat_{<\theta}$ is a full subcategory of $\Cat_{<\theta+1}$,
        \item For a limit ordinal $\lambda$,
            $$
            \Cat_{<\lambda} \simeq \varcolim_{\theta<\lambda}\Cat_{<\theta}
            $$
    \end{itemize}
\end{lemma}
\begin{proof}
    That $\Cat_{<0}\simeq\{*\}$ and $\Cat_{<\theta+1} \simeq (\Cat_{<\theta})\Cat$ follow by definition.
    For the limit case, suppose by transfinite induction that $\Cat_{<\theta}\subseteq\Cat_{<\theta'}$ for all $\theta<\theta'<\lambda$.
    Then,
    \[
        \varcolim_{\theta<\lambda}\Cat_{<\theta} \simeq \bigcup_{\theta<\lambda}\Cat_{<\theta} = \Cat_{<\lambda}
    \]
    as desired.
\end{proof}
\begin{lemma}\label{lem:cat-rank-distinct}
    For every ordinal $\theta$, there is an $(\infty, \infty)$-category $\sC$ such that $\rank\sC = \theta$; that is, $\rank\sC < \theta + 1$ but $\rank\sC\not<\theta$.
\end{lemma}
\begin{proof}
    We prove this by transfinite induction.
    For $\theta = 0$, we take $\sC = \varnothing$.
    Indeed, $\rank\sC < 1$ is vacuous, and $\rank\sC\not<0$ because $\sC\not\simeq*$.
    
    Suppose we have an $(\infty, \infty)$-category $\sD$ such that $\rank\sD = \theta$.
    Then, $\rank\sC = \theta + 1$ for $\sC := \Sigma\sD$.
    
    Finally, suppose $\lambda$ is a limit ordinal such that for every $\theta < \lambda$, there exists an $(\infty, \infty)$-category $\sD^\theta$ such that $\rank\sD^\theta = \theta$.
    Then, take $\sC := \coprod_{\theta<\lambda}\sD^\theta$.
    
    Let $x, y\in\sC$.
    If $x\in\sD^\theta$ and $y\in\sD^{\theta'}$ with $\theta\neq\theta'$, then $\rank\Hom_\sC(x, y) = \rank\varnothing = 0 < \lambda$.
    Otherwise, $\rank\Hom_\sC(x, y) = \rank\Hom_{\sD^\theta}(x, y) < \lambda$.
    In particular, $\rank\sC < \lambda + 1$.
    On the other hand, $\rank\sC\not<\theta$ for all $\theta < \lambda$ since $\sD^\theta$ is a (full) subcategory of $\sC$, and $\rank\sD^\theta\not<\theta$.
    Therefore, $\rank\sC \not< \lambda$, proving that $\rank\sC = \lambda$, as desired.
\end{proof}

\begin{proposition}\label{prop:enr-init-dnt}
    Ad\'amek's construction of an initial algebra for $(-)\Cat$ over the category $\SymMon_\infty$ does not terminate.
\end{proposition}
\begin{proof}
    The $\theta$th stage of Ad\'amek's construction yields $\Cat_{<\theta}$ by \cref{lem:cat-rank}.
    Therefore, the proposition follows from \cref{lem:cat-rank-distinct}.
\end{proof}

The failure of Ad\'amek's construction to terminate is purely a size issue.
For instance, let $(-)\Cat^{<\omega}$ denote the subfunctor of $(-)\Cat$ that sends $\sV$ to the full subcategory $\sV\Cat^{<\omega}$ of $\sV\Cat$ spanned by those $\sV$-enriched categories with finitely many equivalence classes of objects (that is, the underlying space of objects has finitely many path-connected components).
Then, Ad\'amek's construction for $(-)\Cat^{<\omega}$ terminates after $\omega$ steps, and the initial algebra consists of those finite-dimensional higher categories with finitely many equivalence classes of $k$-morphisms for each $k\geq0$.

This phenomenon can be shown more generally:

\begin{lemma}\label{lem:lambda-small-noeth}
    Fix a regular cardinal $\lambda$.
    Let $\sC$ be an $(\infty, \infty)$-category such that the set of equivalence classes of objects of $\sC$ is $\lambda$-small, and $\rank\Hom_\sC(x, y) < \lambda$ for all $x, y\in\sC$.
    Then, $\rank\sC < \lambda$. 
\end{lemma}
\begin{proof}
    For $x, y\in\sC$, let $\theta_{x,y} < \lambda$ such that $\rank\Hom_\sC(x, y) < \theta_{x, y}$; such an ordinal exists because a regular cardinal is necessarily a limit ordinal.
    Then, let $\theta := \sup_{x,y\in\sC}\theta_{x, y}$.
    Note that if $x\simeq x'$ and $y\simeq y'$, then $\theta_{x,y} = \theta_{x',y'}$.
    Since $\sC$ has fewer than $\lambda$ objects up to equivalence, it follows from the fact that $\lambda$ is a regular cardinal that $\theta < \lambda$, and therefore also that $\theta + 1 < \lambda$.
    Therefore, $\rank\sC < \theta + 1 < \lambda$, as desired.
\end{proof}

\begin{proposition}\label{prop:enr-init-smallterminate}
    For a regular cardinal $\lambda$, let $(-)\Cat^{<\lambda}$ denote the subfunctor of $(-)\Cat : \SymMon_\infty \to \SymMon_\infty$ that associates to a symmetric monoidal category $\sV$ the full subcategory $\sV\Cat^{<\lambda}$ of $\sV\Cat$ spanned by those $\sV$-enriched categories such that the set of path-connected components of its underlying space of objects is $\lambda$-small.
    Then, Ad\'amek's construction of an initial algebra for $(-)\Cat^{<\lambda}$ over $\SymMon_\infty$ terminates after no fewer than $\lambda$ steps.
\end{proposition}
\begin{proof}
    For an ordinal $\theta$, let $\Cat_{<\theta}^{<\lambda}$ denote the full subcategory of $\Cat_{<\theta}$ on those $(\infty, \infty)$-categories $\sC$ such that the set of equivalence classes of $k$-morphisms is $\lambda$-small for every $k\geq0$.
    Then, $\Cat_{<\theta}^{<\lambda}$ can be constructed by transfinite induction, analogous to \cref{lem:cat-rank}:
    \begin{itemize}
        \item $\Cat_{<0}^{<\lambda} \simeq \Grpd_{-2} \simeq \{*\}$, which is the initial object in $\SymMon_\infty$,
        \item $\Cat_{<\theta + 1}^{<\lambda} \simeq (\Cat_{<\theta}^{<\lambda})\Cat^{<\lambda}$,
        \item For a limit ordinal $\mu$,
            $$
            \Cat_{<\mu}^{<\lambda} \simeq \varcolim_{\theta<\mu}\Cat_{<\theta}^{<\lambda}
            $$
    \end{itemize}
    Following the proof of \cref{lem:cat-rank-distinct}, there still exists $\sC\in\Cat_{<\theta+1}^{<\lambda}$ such that $\sC\notin\Cat_{<\theta}^{<\lambda}$, so long as $\theta < \lambda$.
    However, \cref{lem:lambda-small-noeth} shows that $\Cat_{<\lambda}^{<\lambda}\subseteq\Cat_{<\theta}^{<\lambda}$ is an equivalence for all $\theta > \lambda$.
    
    Therefore, Ad\'amek's construction terminates in exactly $\lambda$ steps, as desired, and $\Cat_{<\lambda}^{<\lambda}$ carries the structure of an initial algebra for $(-)\Cat^{<\lambda}$ over $\SymMon_\infty$.
\end{proof}

We conclude this subsection by proving that $(-)\Cat$ has an initial algerba over $\SymMon_\infty$.

\begin{definition}\label{def:Noeth}
    A \emph{parallel morphism tower} $(\vec\alpha,\vec\beta)$ in an $(\infty, \infty)$-category $\sC$ is a (countable) sequence of pairs
    $$
    (\alpha_0, \beta_0), (\alpha_1, \beta_1), (\alpha_2, \beta_2), \dots
    $$
    where $\alpha_0,\beta_0$ are objects of $\sC$, and $\alpha_{n+1}$ and $\beta_{n+1}$ are parallel $(n+1)$-morphisms $\alpha_n\to\beta_n$ in $\sC$ for all $n\geq0$.
    
    Say that an $(\infty, \infty)$-category $\sC$ is \emph{Noetherian} if for any parallel morphism tower $(\vec\alpha,\vec\beta)$, there exists $N\gg0$ such that $\Hom_\sC(\alpha_N, \beta_N) \simeq *$.
    
    Denote by $\Cat_{(\infty, \infty)}^{\mathrm{Noeth}}$ the full subcategory of $\Cat_{(\infty, \infty)}$ spanned by the Noetherian $(\infty, \infty)$-categories.
\end{definition}

\begin{lemma}\label{lem:Noeth}
    For an $(\infty, \infty)$-category $\sC$, the following are equivalent:
    \begin{enumerate}[label={(\roman*)}]
        \item\label{it:Noeth}
            $\sC$ is Noetherian,
        \item\label{it:loc-Noeth}
            $\sC$ is locally Noetherian, in the sense that $\Hom_\sC(x, y)$ is Noetherian for all $x,y\in\sC$,
        \item\label{it:small-rank}
            $\sC$ has small rank, in that $\rank\sC < \theta$ for some ordinal $\theta$,
        \item\label{it:loc-small-rank}
            $\sC$ locally has small rank, in that $\Hom_\sC(x, y)$ has small rank for all $x, y\in\sC$.
    \end{enumerate}
\end{lemma}
\begin{proof}
    The equivalence between \ref{it:Noeth} and \ref{it:loc-Noeth} follows by definition.
    
    Note that \ref{it:small-rank} certainly implies \ref{it:loc-small-rank}: if $\rank\sC < \theta$, then $\rank\sC < \theta+1$, and therefore $\rank\Hom_\sC(x, y) < \theta$ for all $x, y\in\sC$.
    Conversely, if for all $x, y\in\sC$ there exists an ordinal $\theta_{x,y}\gg0$ such that $\rank\Hom_\sC(x, y) < \theta_{x, y}$, choose $\lambda\gg0$ such that the set of equivalence classes of objects in $\sC$ is $\lambda$-small, and such that $\lambda \geq\theta_{x, y}$ for all $x, y\in\sC$.
    Then, $\rank\sC < \lambda$ by \cref{lem:lambda-small-noeth}.
    This proves that \ref{it:small-rank} is equivalent to \ref{it:loc-small-rank}.
    
    Since the singleton $*$ is certainly Noetherian, and locally Noetherian $(\infty, \infty)$-categories are Noetherian, it follows by transfinite induction on the rank that every $(\infty, \infty)$-category $\sC$ with small rank is Noetherian.
    This shows that \ref{it:small-rank} implies \ref{it:Noeth}.
    
    To prove the converse, suppose $\sC$ does not have small rank.
    Then, $\sC$ does not locally have small rank, so there must exist $\alpha_0, \beta_0\in\sC$ such that $\Hom_\sC(\alpha_0, \beta_0)$ does not have small rank.
    Proceeding recursively, we obtain a parallel morphism tower $(\vec\alpha,\vec\beta)$ where each $\Hom_\sC(\alpha_n, \beta_n)$ does not have small rank.
    In particular, $\Hom_\sC(\alpha_n, \beta_n)\not\simeq*$ for every $n\geq0$.
    Therefore, if $\sC$ does not have small rank, then $\sC$ is not Noetherian, completing the proof.
\end{proof}

\begin{theorem}\label{thm:enr-init-bigterminate}
    $\Cat_{(\infty, \infty)}^{\mathrm{Noeth}}$ carries the structure of an initial algebra for $(-)\Cat$ over $\SymMon_\infty$.
\end{theorem}
\begin{proof}
    By \cref{lem:Noeth}, the canonical inclusion $\Cat_{(\infty, \infty)}^{\mathrm{Noeth}}\subseteq(\Cat_{(\infty, \infty)}^{\mathrm{Noeth}})\Cat$ is an equivalence.

    By expanding universes, let $\Lambda$ denote the large ordinal of all (small) ordinals.
    Then, \cref{lem:Noeth} implies that $\Cat_{(\infty, \infty)}^{\mathrm{Noeth}}$ is the $\Lambda$-filtered colimit
    $$
    \Cat_{(\infty, \infty)}^{\mathrm{Noeth}} = \bigcup_\theta\Cat_{<\theta} \simeq \varcolim_{\theta<\Lambda}\Cat_{<\theta}
    $$
    in $\SymMon_\infty$, which by \cref{lem:cat-rank} is precisely $\Lambda$ stages of Ad\'amek's initial algebra construction.
    Since the construction terminates after $\Lambda$ steps by the previous discussion, the theorem follows from \cref{thm:adamek}.
\end{proof}
\begin{remark}\label{rem:enr-init-bigterminate}
    Although Ad\'amek's construction in this case requires a large colimit, this colimit is small relative to an expanded universe, and \cref{thm:adamek} applies also to categories that are small (relative to the universe of discourse).
\end{remark}

%% file: src/appendix.tex
%%!TEX root = ../article.tex

\section{Technicalities on lax algebras}\label{adx}

\input{src/appendix/propagation}

\input{src/appendix/limits}

%% file: src/appendix/propagation.tex
%%!TEX root = ../../article.tex

\subsection{Constructing the propagation endofunctor on lax algebras}\label{adx:prop}

This appendix is dedicated to formally constructing the propagation endofunctor on $\sK^\lax(F)$ described in \cref{def:prop}.
To do so, we introduce an auxiliary category:

\begin{definition}\label{def:SKF}
    Define the category $\sK^\square(F)$ to be the (homotopy) pullback
    $$
    \begin{tikzcd}
        \sK^\square(F) \ar[r]\ar[d]\ar[dr, phantom, "\lrcorner" very near start]
        & \sK^{\Lambda^0[2]^\triangleright} \ar[d, two heads] \\
        \sK \ar[r, "{(F,\id)}"']
        & \sK\times\sK
    \end{tikzcd}
    $$
\end{definition}

Note that $\Lambda^0[2]^\triangleright\cong\Delta[1]\times\Delta[1]$ is the walking commutative square, so $\sK^\square(F)$ consists of commutative squares of the form
$$
\begin{tikzcd}
    E \ar[r]\ar[d]
    & FB \ar[d] \\
    B \ar[r]
    & C
\end{tikzcd}
$$
There is an evident forgetful functor $U:\sK^\square(F)\to\sK^\lax(F)$ given on objects by the mapping
$$
\left\{\begin{tikzcd}
    E \ar[r]\ar[d]
    & FB \ar[d] \\
    B \ar[r]
    & C
\end{tikzcd}\right\}
\qquad \mapsto \qquad
\left\{\begin{tikzcd}
    FB \\
    E \ar[u]\ar[d] \\
    B
\end{tikzcd}
\right\}
$$
We can construct a precursor to the propagation endofunctor on $\sK^\lax(F)$ through $\sK^\square(F)$, giving a functor $\Pi^\square:\sK^\square(F)\to\sK^\lax(F)$ and a natural transformation $\eta^\square:U\Rightarrow\Pi^\square$.
Intuitively, this functor and transformation come from the mapping
\begin{align}
    \left\{\begin{tikzcd}[ampersand replacement=\&]
        E \ar[r, "r"]\ar[d]\ar[dr, gray]
        \& FB \ar[d] \\
        B \ar[r, "c"']
        \& C
    \end{tikzcd}\right\}
    \qquad \mapsto \qquad
    \left\{\begin{tikzcd}[ampersand replacement=\&]
        FB \ar[r, "Fc"]
        \& FC \\
        E \ar[u, "r"]\ar[r, "r"]\ar[d]\ar[ur, gray]\ar[dr, gray]
        \& FB \ar[u, "Fc"']\ar[d] \\
        B \ar[r, "c"']
        \& C
    \end{tikzcd}\right\}
    \label{eq:pfi}
\end{align}
To formalise this construction, consider the maps $\bar r,\bar c:\Delta[1]\to\Lambda^0[2]^\triangleright$ which classify the upper edge (corresponding to $E\to FB$ in the diagram) and lower edge (corresponding to $c:B\to C$ in the diagram), respectively, then we can define a functor $\sK^\square(F)\to\sK^{\Lambda^1[2]}$ via
$$
\begin{tikzcd}
    \sK^\square(F) \ar[drr, bend left, "{F\bar c^*}"]\ar[ddr, bend right, "{\bar r^*}"']\ar[dr, dashed, "\exists!"] \\[-1em]
    &[-1em] \sK^{\Lambda^1[2]} \ar[r]\ar[d]\ar[dr, phantom, "\lrcorner" very near start]
    & \sK^{\Delta[1]} \ar[d] \\
    & \sK^{\Delta[1]} \ar[r]
    & \sK
\end{tikzcd}
$$
Explicitly, the dashed arrow describes the mapping
$$
\left\{\begin{tikzcd}
    E \ar[r, "r"]\ar[d]\ar[dr, gray]
    & FB \ar[d] \\
    B \ar[r, "c"']
    & C
\end{tikzcd}\right\}
\qquad\mapsto\qquad
\left\{\begin{tikzcd}
    & FC \\
    E \ar[r, "r"]
    & FB \ar[u, "{Fc}"']
\end{tikzcd}\right\}
$$
However, this mapping does \emph{not} define a composite for this sequence of arrows.
Since $\Lambda^1[2]\hookrightarrow\Delta[2]$ is inner anodyne, and $\sK$ is a quasicategory, the map $\sK^{\Delta[2]}\to\sK^{\Lambda^1[2]}$ is a trivial inner fibration.
Therefore, we can find a section $\sK^{\Lambda^1[2]}\to\sK^{\Delta[2]}$ by solving the lifting problem
$$
\begin{tikzcd}
    \emptyset \ar[r]\ar[d, tail]
    & \sK^{\Delta[2]} \ar[d, two heads, "\sim" sloped] \\
    \sK^{\Lambda^1[2]} \ar[ur, dashed, "\exists"]\ar[r, equal]
    & \sK^{\Lambda^1[2]}
\end{tikzcd}
$$
This provides a functorial choice of composites to the diagram above, and particular provides a functor $\sK^\square(F)\to\sK^{\Delta[2]}$.

In particular, since $\Delta[1]\times\Delta[1]\cong\Lambda^0[2]^\triangleright = (\Delta[1]^\triangleright)\sqcup_{\Delta[0]^\triangleright}(\Delta[1]^\triangleright) = \Delta[2]\sqcup_{\Delta[1]}\Delta[2]$ is obtained by gluing two triangles along their hypotenuse, we can extend the above section to define a map $\sK^\square(F)\to\sK^{\Delta[1]\times\Delta[1]}$ corresponding to the mapping
$$
\left\{\begin{tikzcd}
    E \ar[r, "r"]\ar[d]\ar[dr, gray]
    & FB \ar[d] \\
    B \ar[r, "c"']
    & C
\end{tikzcd}\right\}
\qquad\mapsto\qquad
\left\{\begin{tikzcd}
    FB \ar[r, "{Fc}"]
    & FC \\
    E \ar[u, "r"]\ar[r, "r"]\ar[ur, gray]
    & FB \ar[u, "{Fc}"']
\end{tikzcd}\right\}
$$
which is precisely the upper square in the mapping sketched in (\ref{eq:pfi}).
By gluing this square with the forgetful functor $\sK^\square(F)\to\sK^{\Lambda^0[2]^\triangleright}=\sK^{\Delta[1]\times\Delta[1]}$, which describes the lower square in (\ref{eq:pfi}), we obtain a functor
\[
    \sK^\square(F)\to\sK^{\Delta[1]\times\Delta[1]}\times_{\sK^{\Delta[1]}}\sK^{\Delta[1]\times\Delta[1]} \cong \sK^{\Lambda^0[2]\times\Delta[1]} = \Fun(\Delta[1],\sK^{\Lambda^0[2]})
\]
describing precisely the mapping sketched in (\ref{eq:pfi}).
Moreover, it follows from the construction that the adjunct $\sK^\square(F)\times\Delta[1]\to\sK^{\Lambda^0[2]}$ factors through the forgetful functor $\sK^\lax(F)\to\sK^{\Lambda^0[2]}$.

The resulting functor $\sK^\square(F)\times\Delta[1]\to\sK^\lax(F)$ corresponds to a map $\eta^\square:\Delta[1]\to\Fun(\sK^\square(F),\sK^\lax(F))$.
In particular, $\eta^\square$ classifies a natural transformation between functors $\sK^\square(F)\to\sK^\lax(F)$ whose domain, by construction, is precisely the forgetful functor $U$.
\begin{definition}\label{def:spu}
    Let $\Pi^\square$ be the codomain of the natural transformation constructed above, and denote the natural transformation itself by $\eta^\square: U\Rightarrow\Pi^\square$.
\end{definition}

We can now use the above construction to create the propagation endofunctor on $\sK^\lax(F)$, as well as its unit.
\begin{definition}\label{def:p+u}
    Let $\sK$ be finitely cocomplete.
    By \cite[Proposition 4.2.2.7]{lurie}, taking colimits defines a functor $\varcolim:\sK^{\Lambda^0[2]}\to\sK^{\Lambda^0[2]^\triangleright}$.
    It restricts to a functor $\varcolim:\sK^\lax(F)\to\sK^\square(F)$, which is a section of the forgetful functor $U:\sK^\square(F)\to\sK^\lax(F)$.

    Define the \emph{propagation endofunctor} $\Pi$ on $\sK^\lax(F)$ to be the composite
    \[
        \sK^\lax(F)\xrightarrow{\varcolim}\sK^\square(F)\xrightarrow{\Pi^\square}\sK^\lax(F)
    \]
    This functor then admits a \emph{unit} given by the composite
    \[
        \eta:\Delta[1]\xrightarrow{\eta^\square}\Fun(\sK^\square(F),\sK^\lax(F))\xrightarrow{\varcolim^*}\Fun(\sK^\lax(F),\sK^\lax(F))
    \]
    Indeed, this classifies a natural transformation $\Id\Rightarrow\Pi$ because $\varcolim$ is a section of $U$.
\end{definition}

%% file: src/appendix/limits.tex
%%!TEX root = ../../article.tex

\subsection{Colimits of lax algebras}\label{adx:limit}

This appendix is dedicated to computing colimits of lax algebras.
Recall from \cref{rem:tower} that $\sK^\lax(F)$ fits into the pullback square
$$
\begin{tikzcd}
    \sK^\lax(F) \ar[r]\ar[d, "u"']\ar[dr, phantom, "\lrcorner" very near start]
    & \sK^{\Lambda^0[2]} \ar[d] \\
    \sK^{\Delta[1]} \ar[r, "{(F(2),\id)}"']
    & \sK\times\sK^{\Delta[1]}
\end{tikzcd}
$$
The main result of this appendix is the following, which is the key technical result of the paper.

\begin{proposition}\label{prop:laxcol}
    The forgetful functor $\sK^\lax(F)\xrightarrow u\sK^{\Delta[1]}$ reflects colimits.

    Specifically, let $p:J\to\sK^\lax(F)$ be a map of simplicial sets, and say that the lax $F$-algebra at $p_j$ is given by $FB_j\gets E_j\to B_j$.
    Suppose $up:J\to\sK^{\Delta[1]}$ admits a colimit $E_\infty\to B_\infty$.
    Then, we have a cocone of maps $E_j\to FB_j\to FB_\infty$, and so by the universal property of $E_\infty=\varcolim_jE_j$, there is an essentially unique map $E_\infty\to FB_\infty$.
    The resulting lax $F$-algebra $FB_\infty\gets E_\infty\to B_\infty$ is then a colimit of $p:J\to\sK^\lax(F)$.
\end{proposition}

In order to prove this, we will rely on the following technical results:

\begin{lemma}\label{lem:pbc}
    Let $p:J\to\sA\times_\sC\sB$ be a map of simplicial sets into a strict fibre product of quasicategories, and suppose that the composite $\pi_\sA p:J\to\sA$ admits a colimit $\bar\pi:J^\triangleright\to\sA$.
    Then, $p$ admits a colimit in $\sA\times_\sC\sB$ if and only if we can always solve the lifting problem
    \begin{align}
    \begin{tikzcd}[ampersand replacement=\&]
        J\star T \ar[rr]\ar[d]
        \&\& \sB \ar[d] \\
        J^\triangleright\star T\ar[urr, dashed]\ar[r]
        \& \sA \ar[r]
        \& \sC
    \end{tikzcd}
        \label{eq:pbld}
    \end{align}
    where $T$ is any simplicial set, $J\star T\to\sB$ extends $\pi_\sB p:J\to\sB$, and $J^\triangleright\star T\to\sA$ extends $\bar\pi$.
\end{lemma}
\begin{proof}
    If $p$ admits a colimit, then certainly every such lifting problem (\ref{eq:pbld}) can be solved.
    Conversely, suppose every lifting problem (\ref{eq:pbld}) can be solved.
    By taking $T=\emptyset$, we obtain a cocone $\bar p:J^\triangleright\to\sA\times_\sC\sB$ extending $p$.
    The goal is to show that $\bar p$ is indeed a colimit cocone for $p$.
    Therefore, we need to find a dotted arrow fitting in any diagram
    $$
    \begin{tikzcd}
        J\star S \ar[r]\ar[d, hook]
        & J\star T \ar[dr]\ar[d] \\
        J^\triangleright\star S \ar[r]\ar[rr, bend right]
        & J^\triangleright\star T\ar[r, dotted]\ar[dr, dashed]
        & \sA\times_\sC\sB \ar[r]\ar[d]\ar[dr, phantom, "\lrcorner" very near start]
        & \sB \ar[d] \\
        && \sA \ar[r]
        & \sC
    \end{tikzcd}
    $$
    where the map $J\star T\to\sA\times_\sC\sB$ extends $p$, and $J^\triangleright\star S\to\sA\times_\sC\sB$ extends $\bar p$.
    Since we have a colimit cocone $\bar\pi$ in $\sA$, we can find an arrow $J^\triangleright\star T\to\sA$ fitting as the dashed arrow in the above diagram.
    This reduces the problem for finding a dotted arrow into solving a lifting problem (\ref{eq:pbld}), which can be done by assumption.
\end{proof}

\begin{lemma}\label{lem:jx1}
    For all simplicial sets $S, T$, the diagram
    $$
    \begin{tikzcd}
        (S\times\{0\})\star(T\times\{1\}) \ar[r]\ar[d]\ar[dr, phantom, "\ulcorner" very near end]
        & (S\times\Delta[1])\star(T\times\{1\}) \ar[d] \\
        (S\times\{0\})\star(T\times\Delta[1]) \ar[r]
        & (S\star T)\times\Delta[1]
    \end{tikzcd}
    $$
    is a pushout square, which is moreover a homotopy pushout square as the maps are cofibrations.
\end{lemma}
\begin{proof}
    On $n$-cells, the square is given by
    $$
    \resizebox{\textwidth}{!}{\begin{tikzcd}[cells={font=\everymath\expandafter{\the\everymath\displaystyle}}, ampersand replacement=\&]
        (S_n\times\{0\})\sqcup(T_n\times\{1\})\sqcup\coprod_{i+j=n-1}(S_i\times\{0\})\times(T_j\times\{1\}) \ar[r]\ar[d]
        \& (S_n\times\Delta[1]_n)\sqcup(T_n\times\{1\})\sqcup\coprod_{i+j=n-1}(S_i\times\Delta[1]_i)\times(T_j\times\{1\}) \ar[d] \\
        (S_n\times\{0\})\sqcup(T_n\times\Delta[1]_n)\sqcup\coprod_{i+j=n-1}(S_i\times\{0\})\times(T_j\times\Delta[1]_j) \ar[r]
        \& (S_n\times\Delta[1]_n)\sqcup(T_n\times\Delta[1]_n)\sqcup\coprod_{i+j=n-1}(S_i\times T_j\times\Delta[1]_n)
    \end{tikzcd}}
    $$
    Since colimits commute with colimits, it suffices to show that the diagram restricted to each set of coproduct summands forms a pushout square.
    In other words, it suffices to show for all $i+j=n-1$ (where $i,j\geq-1$ and we take $K_{-1} := *$ for any simplicial set $K$) that the square
    $$
    \begin{tikzcd}
        (S_i\times\{0\})\times(T_j\times\{1\}) \ar[r]\ar[d]
        & (S_i\times\Delta[1]_i)\times(T_j\times\{1\}) \ar[d] \\
        (S_i\times\{0\})\times(T_j\times\Delta[1]_j) \ar[r]
        & S_i\times T_j\times\Delta[1]_{i+j+1}
    \end{tikzcd}
    $$
    is a pushout square.
    This is trivial if $i=-1$ or $j=-1$, so suppose $i,j\geq0$.
    Since $\sSet$ is cartesian closed, products commute with colimits, which allows us to reduce further to showing that
    $$
    \begin{tikzcd}
        * \ar[r]\ar[d]
        & \Hom_{\bDelta}([i],[1]) \ar[d] \\
        \Hom_{\bDelta}([j],[1]) \ar[r]
        & \Hom_{\bDelta}([i+j+1],[1])
    \end{tikzcd}
    $$
    is a pushout square.
    Note that maps $[i]\to[1]$ correspond to integers $0\leq c\leq i+1$, where $c$ indicates the first index of the map that is sent to $1$ (and $c=i+1$ means that the map is constant at zero).
    With this interpretation, the top map picks out the morphism $[i]\to[1]$ corresponding to the integer $c=0$, while the vertical map on the left picks out the zero map $[j]\to[1]$.
    Observing that we have a pushout square amounts to observing that a morphism $[i+j+1]\to[1]$ falls into one of the following three cases:
    \begin{itemize}
        \item it corresponds to a cut $0\leq c<j+1$, in which case it comes from a nonzero morphism $[j]\to[1]$
        \item it corresponds to a cut $j+1<c\leq i+j+2$, in which case it comes from a morphism $[i]\to[1]$ that starts at zero (the vertical map on the right shifts the index of the cut up by $j+1$)
        \item it corresponds to the cut $c=j+1$, in which case it simultaneously comes from the constant zero morphism $[j]\to[1]$ and the constant one morphism $[i]\to[1]$
    \end{itemize}
\end{proof}

We moreover recall the following result:
\begin{lemma}\cite[Lemma 2.1.2.3]{lurie}\label{lem:psj}
    Let $A_0\subseteq A$ and $B_0\subseteq B$ be inclusions such that either $A_0\subseteq A$ is right anodyne, or $B_0\subseteq B$ is left anodyne.
    Then, the inclusion
    \[
        (A_0\star B)\sqcup_{A_0\star B_0}(A\star B_0) \hookrightarrow A\star B
    \]
    is inner anodyne.
\end{lemma}

\begin{proof}[Proof of \cref{prop:laxcol}]
    Suppose we have a diagram $p:J\to\sK^\lax(F)$ such that $up:J\to\sK^{\Delta[1]}$ admits a colimit $\bar p_u:J^\triangleright\to\sK^{\Delta[1]}$.
    By \cref{lem:pbc}, it suffices to show for any simplicial set $T$ that we can find a lift for any problem
    \begin{align}
        \begin{tikzcd}[ampersand replacement=\&]
            J\star T \ar[d, hook]\ar[r, "q"]
            \& \sK^\lax(F) \ar[r]
            \& \sK^{\{1\gets0\to2\}} \ar[d]\ar[r]\ar[dr, phantom, "\lrcorner" very near start]
            \& \sK^{\{1\gets0\}} \ar[d] \\
            J^\triangleright\star T \ar[r]\ar[urr, dashed]
            \& \sK^{\{0\to2\}} \ar[r, "{(F(2),\id)}"']
            \& \sK^{\{1\}}\times\sK^{\{0\to2\}} \ar[r]
            \& \sK^{\{1\}}\times\sK^{\{0\}}
        \end{tikzcd}
        \label{eq:rlpc}
    \end{align}
    where $q$ extends the composite $J\xrightarrow p\sK^\lax(F)$, and the leftmost arrow on the bottom extends $\bar p_u:J^\triangleright\to\sK^{\Delta[1]}$.
    Since we have the pullback square on the right, it suffices to find a lift $J^\triangleright\star T\to\sK^{\{0\to1\}}$ to the upper right corner.
    By currying, we are finding a suitable map $(J^\triangleright\star T)\times\Delta[1]\to\sK$.
    
    For the sake of clarity, we will refer to maps of simplicial sets based on an intuitive diagram that they reflect.
    For this purpose, we will denote the lax $F$-algebra $q_j$ at $j\in J\subseteq J\star T$ by $FB_j\gets E_j\to B_j$, and we will denote the lax $F$-algebra $q_t$ at $t\in T\subseteq J\star T$ by $FC_t\gets D_t\to C_t$.
    Similarly, denote the colimit $p_u(\infty)\in\sK^{\{0\to2\}}$ by $E_\infty\to B_\infty$.
    Then, the desired lift $(J^\triangleright\star T)\times\Delta[1]\to\sK$ reflects the diagram
    \begin{align}
        \left\{\begin{tikzcd}[ampersand replacement=\&, sep=small]
            FB_j \ar[rr]\ar[dr]
            \&\& FC_t \\
            \& FB_\infty \ar[ur] \\
            E_j \ar[uu]\ar[rr]\ar[dr]
            \&\& D_t \ar[uu] \\
            \& E_\infty \ar[uu, crossing over]\ar[ur]
        \end{tikzcd}\right\}
        \label{eq:lft0}
    \end{align}
    The first goal is to produce a map reflecting the diagram
    \begin{align}
        \left\{\begin{tikzcd}[sep=small, ampersand replacement=\&]
            \& FB_\infty \ar[dr] \\
            E_j \ar[ur]\ar[dr]\ar[r, dashed]
            \& E_\infty \ar[u, dashed]\ar[d, dashed]\ar[r, dashed]
            \& FC_t \\
            \& D_t \ar[ur]
        \end{tikzcd}\right\}
        \label{eq:goal}
    \end{align}
    from which the dashed arrows can be recovered by the universal property of $E_\infty=\varcolim_jE_j$ (note that the dashed arrows $E_\infty\to D_t$ are already provided from the bottom row of (\ref{eq:rlpc}) via the map $\{\infty\}\star T\subseteq J^\triangleright\star T\to\sK^{\{0\to2\}}\to\sK^{\{0\}}$).

    The commutative diagram
    \begin{align}
        \left\{\begin{tikzcd}[sep=small, ampersand replacement=\&]
            FB_j \ar[r]
            \& FC_t \\
            E_j \ar[u]\ar[r]\ar[ur, gray]
            \& D_t \ar[u]
        \end{tikzcd}\right\}
        \label{eq:lft1}
    \end{align}
    is obtained by the map $(J\star T)\times\Delta[1]\to\sK$ given as the adjunct of the top row of (\ref{eq:rlpc}).
    By \cref{lem:jx1}, we can write
    \[
        (J\star T)\times\Delta[1] &= \big((J\times\{0\})\star(T\times\Delta[1])\big)\sqcup_{(J\times\{0\})\star(T\times\{1\})}\big((J\times\Delta[1])\star(T\times\{1\})\big)
    \]
    In particular, we can isolate the upper-left triangle of (\ref{eq:lft1}) by restricting to the simplicial subset $(J\times\Delta[1])\star(T\times\{1\})\hookrightarrow(J\star T)\times\Delta[1]\to\sK$.
    We also have a map $J\star\{\infty\}\star T\to\sK$ obtained from the bottom row of (\ref{eq:rlpc}) as the composite $J^\triangleright\star T\to\sK^{\{1\}}\times\sK^{\{0\to2\}}\to\sK^{\{1\}}$, which reflects the diagram
    $$
    \left\{\begin{tikzcd}[sep=small]
        & FB_\infty \ar[dr] \\
        FB_j \ar[ur]\ar[rr]
        && FC_t
    \end{tikzcd}\right\}
    $$
    Gluing with the upper left triangle of (\ref{eq:lft1}), we can produce a map $(J\times\Delta[1])\star\{\infty\}\star(T\times\{1\})\to\sK$ reflecting the diagram
    \begin{align}
        \left\{\begin{tikzcd}[ampersand replacement=\&, sep=small]
            FB_j \ar[r]\ar[drr, bend right=10]
            \& FB_\infty \ar[dr] \\
            E_j \ar[u]\ar[ur]\ar[rr]
            \&\& FC_t
        \end{tikzcd}\right\}
        \label{eq:lft2}
    \end{align}
    as a solution to the lifting problem
    $$
    \begin{tikzcd}
        \big((J\times\Delta[1])\star(T\times\{1\})\big)\sqcup_{(J\times\{1\})\star(T\times\{1\})}\big((J\times\{1\})\star\{\infty\}\star(T\times\{1\})\big) \ar[r]\ar[d, hook]
        & \sK \\
        (J\times\Delta[1])\star\{\infty\}\star(T\times\{1\}) \ar[ur, dashed]
    \end{tikzcd}
    $$
    By applying \cref{lem:psj} to the right anodyne map $J\times\{1\}=J\times\Lambda^1[1]\hookrightarrow J\times\Delta[1]$ and the inclusion $T\times\{1\}\hookrightarrow\{\infty\}\star(T\times\{1\})$, we see that the vertical map is inner anodyne.
    Therefore, since $\sK$ is a quasicategory, it follows that such a lift indeed exists.

    Now, glue the map for (\ref{eq:lft2}) to the bottom-right triangle of (\ref{eq:lft1}) to produce a map
    \[
        \big((J\times\{0\})\star(T\times\Delta[1])\big)\sqcup_{(J\times\{0\})\star(T\times\{1\})}\big((J\times\Delta[1])\star\{\infty\}\star(T\times\{1\})\big)\to\sK
    \]
    In particular, we can restrict this map to the simplicial subset
    \[
        &\big((J\times\{0\})\star(T\times\Delta[1])\big)\sqcup_{(J\times\{0\})\star(T\times\{1\})}\big((J\times\{0\})\star\{\infty\}\star(T\times\{1\})\big) \\
        &\cong(J\times\{0\})\star\Big((T\times\Delta[1])\sqcup_{T\times\{1\}}\big(\{\infty\}\star(T\times\{1\})\big)\Big)
    \]
    reflecting the subdiagram
    $$
    \left\{\begin{tikzcd}[sep=small]
        & FB_\infty \ar[dr] \\
        E_j \ar[ur]\ar[rr]\ar[dr]
        && FC_t \\
        & D_t \ar[ur]
    \end{tikzcd}\right\}
    $$
    which is precisely the perimeter of (\ref{eq:goal}) required to invoke the universal property of $E_\infty$.
    Indeed, since $E_\infty$ is a colimit of the diagram $J\xrightarrow p\sK^\lax(F)\to\sK^{\{1\gets0\to2\}}\to\sK^{\{0\}}$, it follows that we can find a dashed morphism fitting in the diagram
    $$
    \begin{tikzcd}
        (J\times\{0\})\star(T\times\{0\}) \ar[r]\ar[d]
        & (J\times\{0\})\star\Big((T\times\Delta[1])\sqcup_{T\times\{1\}}\big(\{\infty\}\star(T\times\{1\})\big)\Big) \ar[d]\ar[ddr, bend left] \\
        (J^\triangleright\times\{0\})\star(T\times\{0\}) \ar[r]\ar[drr, bend right=10]
        & (J^\triangleright\times\{0\})\star\Big((T\times\Delta[1])\sqcup_{T\times\{1\}}\big(\{\infty\}\star(T\times\{1\})\big)\Big) \ar[dr, dashed] \\
        && \sK
    \end{tikzcd}
    $$
    where the bottom map $(J^\triangleright\times\{0\})\star(T\times\{0\})\to\sK$ is the projection of the bottom row of (\ref{eq:rlpc}) onto $\sK^{\{0\}}$ that describes the complex of morphisms
    $$
    \left\{\begin{tikzcd}[sep=small]
        & E_\infty \ar[dr] \\
        E_j \ar[ur]\ar[rr]
        && D_t
    \end{tikzcd}\right\}
    $$
    This dashed morphism precisely recovers the diagram (\ref{eq:goal}).
    To obtain the desired diagram (\ref{eq:lft0}), we glue this morphism with the map $(J\times\Delta[1])\star\{\infty\}\star(T\times\{1\})\to\sK$ describing (\ref{eq:lft2}).
    Indeed, we get a pushout square
    \begin{align}
        \resizebox{\textwidth}{!}{\begin{tikzcd}[ampersand replacement=\&]
            (J\times\{0\})\star\{\infty\}\star(T\times\{1\}) \ar[r]\ar[d]\ar[dr, phantom, "\ulcorner" very near end]
            \& (J\times\Delta[1])\star\{\infty\}\star(T\times\{1\}) \ar[d]\ar[ddr, bend left] \\
            (J^\triangleright\times\{0\})\star\Big((T\times\Delta[1])\sqcup_{T\times\{1\}}\big(\{\infty\}\star(T\times\{1\})\big)\Big) \ar[r]\ar[drr, bend right=10]
            \& (J^\triangleright\star T)\times\Delta[1] \ar[dr, dashed] \\
            \&\& \sK
        \end{tikzcd}}
        \label{eq:pos}
    \end{align}
    and the dashed arrow precisely reflects the diagram (\ref{eq:lft0}), meaning its adjunct $J^\triangleright\star T\to\sK^{\{0\to1\}}$ pulls back to give precisely a lift in (\ref{eq:rlpc}), as desired.

    To see that (\ref{eq:pos}) is indeed a pushout square, note that by expanding the pushout on the bottom left corner and using that joins preserve pushouts, this is equivalent to showing that the diagram
    $$
    \resizebox{\textwidth}{!}{\begin{tikzcd}[ampersand replacement=\&]
        \& (J\times\{0\})\star\{\infty_1\}\star(T\times\{1\}) \ar[d]\ar[r]
        \& (J\times\Delta[1])\star\{\infty_1\}\star(T\times\{1\}) \ar[dd] \\
        (J^\triangleright\times\{0\})\star(T\times\{1\}) \ar[r]\ar[d]
        \& (J^\triangleright\times\{0\})\star\{\infty_1\}\star(T\times\{1\}) \ar[dr] \\
        (J^\triangleright\times\{0\})\star(T\times\Delta[1]) \ar[rr]
        \&\& (J^\triangleright\star T)\times\Delta[1]
    \end{tikzcd}}
    $$
    is a universal cocone diagram, where I have tacitly replaced $\{\infty\}$ with $\{\infty_1\}$, to indicate that its image in $(J^\triangleright\star T)\times\Delta[1]$ lies in the top cell $(J\star\{\infty_1\}\star T)\times\{1\}$.
    Now, notice that the simplicial set $(J\times\{0\})^\triangleright\star\{\infty_1\}\star(T\times\{1\})$ is precisely the simplicial set $(J\times\{0\})\star(\{\infty\}\times\Delta[1])\star(T\times\{1\})$ using the associativity of the join operation and how $\{\infty_0\}\star\{\infty_1\}\cong\{\infty\}\times\Delta[1]$.
    Therefore, the above diagram is precisely the result of pasting the following two pushout squares:
    $$
    \resizebox{\textwidth}{!}{\begin{tikzcd}[ampersand replacement=\&]
        \& (J\times\{0\})\star\{\infty_1\}\star(T\times\{1\}) \ar[r]\ar[d]\ar[dr, phantom, "\ulcorner" very near end]
        \& (J\times\Delta[1])\star\{\infty_1\}\star(T\times\{1\}) \ar[d] \\
        (J^\triangleright\star\times\{0\})\star(T\times\{1\}) \ar[r]\ar[d]
        \& (J\times\{0\})\star(\{\infty\}\times\Delta[1])\star(T\times\{1\}) \ar[r]\ar[dr, phantom, "\ulcorner" very near end]
        \& (J^\triangleright\times\Delta[1])\star(T\times\{1\}) \ar[d] \\
        (J^\triangleright\times\{0\})\star(T\times\Delta[1]) \ar[rr]
        \&\& (J^\triangleright\star T)\times\Delta[1]
    \end{tikzcd}}
    $$
    Indeed, the bottom square is precisely an instance of \cref{lem:jx1}, and the upper square is the result of applying $(-)\star(T\times\{1\})$ to another instance of \cref{lem:jx1}.
    This proves that (\ref{eq:pos}) is indeed a pushout square diagram, and thus we have our desired lift of (\ref{eq:rlpc}).
\end{proof}

A very similar result holds for colimits of coalgebras:

\begin{proposition}\label{prop:cocol}
    Let $\sK$ be an $\infty$-category, and $F:\sK\to\sK$ an endofunctor.
    Then, the forgetful functor $\sK_\co(F)\to\sK$ sending an $F$-coalgebra to its underlying object reflects colimits.
\end{proposition}
\begin{proof}
    Note that $\sK_\co(F)$ is equivalent to the full subcategory of $\sK^\lax(F)$ spanned by those lax $F$-algebras where the lax action is invertible.
    Since the colimit of equivalences in $\sK^{\Delta[1]}$ is an equivalence, the result follows from \cref{prop:laxcol}.
\end{proof}
\begin{corollary}\label{prop:alglim}
    The forgetful functor $\sK(F)\to\sK$ reflects limits.
\end{corollary}

%% file: article.bbl
\begin{thebibliography}{Lam68}

\bibitem[Ad{\'a}74]{adamek}
Ji{\v{r}}{\'\i} Ad{\'a}mek.
\newblock Free algebras and automata realizations in the language of
  categories.
\newblock {\em Commentationes Mathematicae Universitatis Carolinae},
  15(4):589--602, 1974.
\newblock \url{https://dml.cz/handle/10338.dmlcz/105583}.

\bibitem[Ad{\'a}05]{adamek:modern}
Ji{\v{r}}{\'\i} Ad{\'a}mek.
\newblock Introduction to coalgebra.
\newblock {\em Theory and Applications of Categories}, 14(8):157--199, 2005.
\newblock \href{http://www.tac.mta.ca/tac/volumes/14/8/14-08abs.html}{\UrlFont
  TAC:volumes/14/08}.

\bibitem[Ara15]{ara}
Dimitri Ara.
\newblock Higher quasi-categories vs higher {Rezk} spaces.
\newblock {\em Journal of {$K$}-Theory}, 14(3):701--749, 2015.
\newblock \href{https://arxiv.org/abs/1206.4354v3}{\UrlFont arXiv:1206.4354v3}.
\newblock \href {https://doi.org/10.1017/S1865243315000021}
  {\path{doi:10.1017/S1865243315000021}}.

\bibitem[BSP21]{barwick-schommer-pries}
Clark Barwick and Christopher Schommer-Pries.
\newblock On the unicity of the theory of higher categories.
\newblock {\em Journal of the American Mathematical Society}, 34:1011--1058,
  2021.
\newblock \href{https://arxiv.org/abs/1112.0040v6}{\UrlFont arXiv:1112.0040v6}.
\newblock \href {https://doi.org/10.1090/jams/972}
  {\path{doi:10.1090/jams/972}}.

\bibitem[Che07]{cheng}
Eugenia Cheng.
\newblock An {$\omega$}-category with all duals is an {$\omega$}-groupoid.
\newblock {\em Applied Categorical Structures}, 15:439--453, 2007.
\newblock \href {https://doi.org/10.1007/s10485-007-9081-8}
  {\path{doi:10.1007/s10485-007-9081-8}}.

\bibitem[GH15]{gepner-haugseng}
David Gepner and Rune Haugseng.
\newblock Enriched $\infty$-categories via non-symmetric $\infty$-operads.
\newblock {\em Advances in Mathematics}, 279:575--716, 2015.
\newblock \href{https://arxiv.org/abs/1312.3178v4}{\UrlFont arXiv:1312.3178v4}.
\newblock \href {https://doi.org/10.1016/j.aim.2015.02.007}
  {\path{doi:10.1016/j.aim.2015.02.007}}.

\bibitem[Hau15]{haugseng}
Rune Haugseng.
\newblock Rectification of enriched {$\infty$}-categories.
\newblock {\em Algebraic \& Geometric Topology}, 15(4):1931--1982, 2015.
\newblock \href{https://arxiv.org/abs/1312.3881v4}{\UrlFont arXiv:1312.3881v4}.
\newblock \href {https://doi.org/10.2140/agt.2015.15.1931}
  {\path{doi:10.2140/agt.2015.15.1931}}.

\bibitem[HL23]{henry-loubaton}
Simon Henry and Felix Loubaton.
\newblock An inductive model structure for strictg {$\infty$}-categories, 2023.
\newblock Preprint. \href{https://arxiv.org/abs/2301.11424v1}{\UrlFont
  arXiv:2301.11424v1}.

\bibitem[Kel82]{kelly}
G.M. Kelly.
\newblock {\em Basic Concepts of Enriched Category Theory}, volume~64 of {\em
  LMS Lecture Notes in Mathematics}.
\newblock Cambridge University Press, 1982.
\newblock
  \href{http://www.tac.mta.ca/tac/reprints/articles/10/tr10abs.html}{\UrlFont
  TAC:reprints/10}.

\bibitem[Lam68]{lambek}
Joachim Lambek.
\newblock A fixpoint theorem for complete categories.
\newblock {\em Mathematische Zeitschrift}, 103:151--161, 1968.
\newblock \href {https://doi.org/10.1007/BF01110627}
  {\path{doi:10.1007/BF01110627}}.

\bibitem[Lur09]{lurie}
Jacob Lurie.
\newblock {\em Higher Topos Theory}, volume 170 of {\em the Annals of
  Mathematics Studies}.
\newblock Princeton University Press, 2009.
\newblock \href{https://arxiv.org/abs/math/0608040v4}{\UrlFont
  arXiv:math/0608040v4}.

\bibitem[Lur17]{lurie:algebra}
Jacob Lurie.
\newblock Higher algebra, 2017.
\newblock Preprint. \url{https://people.math.harvard.edu/~lurie/papers/HA.pdf}.

\end{thebibliography}
